\documentclass[a4paper]{article}
\usepackage[utf8]{inputenc}
\usepackage{hyperref}
\usepackage[noblocks]{authblk}
\usepackage{amsmath, amsthm, amssymb}
\usepackage{mathtools}
\usepackage{paralist}
\usepackage{type1cm}
\usepackage{type1ec}						
\usepackage{float}							
\usepackage{graphicx}						
\usepackage{bm}								
\usepackage{enumitem}
\usepackage{subfigure}						
\usepackage{color}
\usepackage{xspace}
\usepackage{siunitx}
\usepackage{booktabs}
\usepackage{textcomp}
\usepackage{lmodern}
\usepackage{epstopdf}

\numberwithin{equation}{section}

\newcommand{\MU}{\mu}
\newcommand{\ooc}{1\!/\!c}

\newcommand{\MULTIBAT}{MULTIBAT\xspace}
\newcommand{\BEST}{BEST\xspace}
\newcommand{\PYMOR}{pyMOR\xspace}
\newcommand{\EI}{POD/EI\xspace}
\newcommand{\MOR}{MOR\xspace}
\newcommand{\Li}{lithium\xspace}
\newcommand{\Liion}{lithium-ion\xspace}
\newcommand{\Liions}{lithium ions\xspace}
\newcommand{\figref}[1]{Fig.~\ref{#1}}

\newcommand{\secref}[1]{Sec.~\ref{#1}}
\newcommand{\chem}[3]{\ensuremath{\mathrm{#1}_\mathrm{#2}^\mathrm{#3}}}
\newcommand{\vart}[3]{\ensuremath{#1_\mathrm{#2}^\mathrm{#3}}}
\newcommand{\var}[2]{\ensuremath{#1_\mathrm{#2}}}
\newcommand{\Span}{\operatorname{span}}

\newcommand{\MORReducedTime}{7m 48s}
\newcommand{\MORReducedTimeMinutes}{8}
\newcommand{\MORDetailedTime}{15h 38m 35s}
\newcommand{\MORDetailedTimeHours}{16}
\newcommand{\MORSpeedup}{120}
\newcommand{\MORSpeedupB}{120.3}
\newcommand{\MOROfflineTime}{13h 43m 12s}
\newcommand{\MORCDIM}{183}
\newcommand{\MORPDIM}{69}
\newcommand{\MORBVDIM}{952}
\newcommand{\MORLAMBDADIM}{1027}

\newcommand{\MORCRelErr}{$4.81\cdot 10^{-4}$}
\newcommand{\MORPRelErr}{$4.50\cdot 10^{-3}$}

\newcommand{\MORCEstMax}{$3.61\cdot 10^{-4}$}
\newcommand{\MORCEffMax}{2.89}

\newcommand{\MORCEffMinRec}{1.08}

\newcommand{\MORPEstMax}{$5.55\cdot 10^{-3}$}
\newcommand{\MORPEffMax}{1.45}

\newcommand{\MORPEffMinRec}{3.46}

\newcommand{\MORREALCDIM}{178}

\newcommand{\MORREALPDIM}{67}

\newcommand{\MORREALBVDIM}{924}
\newcommand{\MORREALLAMBDADIM}{997}

\begin{document}
\title{\MULTIBAT: Unified workflow for fast electrochemical 3D simulations of \Liion cells combining virtual stochastic microstructures, electrochemical degradation models and model order reduction}

\author[a]{Julian Feinauer}
\author[b,c]{Simon Hein}
\author[d]{Stephan Rave}
\author[e]{Sebastian~Schmidt}
\author[a,*]{Daniel Westhoff}
\author[e]{Jochen Zausch}
\author[e,f]{Oleg~Iliev}
\author[b,c,g]{Arnulf Latz}
\author[d]{Mario Ohlberger}
\author[a]{Volker~Schmidt}

\affil[a]{Institute of Stochastics, Ulm University, Germany.}
\affil[b]{Institute of Engineering Thermodynamics, German Aerospace Center (DLR), Stuttgart, Germany.}
\affil[c]{Helmholtz Institute for Electrochemical Energy Storage (HIU), Ulm, Germany.\vspace*{-0.75em}}
\affil[d]{Applied Mathematics, Center for Nonlinear Science \& Center for Multiscale Theory and Computation,
          University of Münster, Germany.}
\affil[e]{Fraunhofer Institute for Industrial Mathematics ITWM, Kaiserslautern, Germany.}
\affil[f]{Institute of Mathematics and Informatics, Bulgarian Academy of Science, Sofia, Bulgaria.}
\affil[g]{Institute of Electrochemistry, Ulm University, Germany.}
\affil[*]{Fax:+49 731 50 26349; Tel: +49 731 50 23617; E-mail: daniel.westhoff@uni-ulm.de}

\date{December 5, 2017}
\renewcommand{\Affilfont}{\scriptsize}
\maketitle

\begin{abstract}
We present a simulation workflow for efficient investigations of the interplay between 3D \Liion electrode microstructures and electrochemical performance, with emphasis on lithium plating. Our approach addresses several challenges. First, the 3D microstructures of porous electrodes are generated by a parametric stochastic model, in order to significantly reduce the necessity of tomographic imaging. Secondly, we integrate a consistent microscopic, 3D spatially-resolved physical model for the electrochemical behavior of the lithium-ion cells taking lithium plating and stripping into account. This highly non-linear mathematical model is solved numerically on the complex 3D microstructures to compute the transient cell behavior. Due to the complexity of the model and the considerable size of realistic microstructures even a single charging cycle of the battery requires several hours computing time. This renders large scale parameter studies extremely time consuming. Hence, we develop a mathematical model order reduction scheme. We demonstrate how these aspects are integrated into one unified workflow, which is a step towards computer aided engineering for the development of more efficient \Liion cells.
\end{abstract}

\section{Introduction}
The ubiquity and importance of rechargeable \Liion batteries lead to the increasing demand for physics-based simulation methods that are able to analyze and predict battery behavior. These methods can not only contribute in improving cell design and operation, but they can also greatly support battery research in its understanding of basic mechanisms, like lithium plating that determine battery life and safety, which is yet not well understood.

The electrochemical simulation of \Liion cells goes back to the work of Newman and his co-workers \cite{Doyle1993,Fuller1994,Newman2003}. Their simulation methodology is based on the porous electrode theory developed by Newman \cite{Newman2012}. This model approach neglects the details of electrode microstructures and describes them as a homogeneous medium where electrolyte and the solid material coexist at every point. The most commonly used model of Newman only considers the through-direction of the battery. It takes into account the diffusion of \Liions into the active material by assuming a spherical, microscopic particle of average size in each discretization point in which a one-dimensional diffusion equation is solved. Hence this model is sometimes called a pseudo-2d (P2D) model \cite{Santhanagopalan2006}. 
There exist many applications for Newman-type models like the study of cell behavior as well as degradation \cite{Santhanagopalan2006,Arora1998,Remmlinger2014,Tippmann2013}. But the main drawback of these models is that the complex electrode microstructures are only approximately accounted for by a few aggregated parameters: the thickness of the electrode $L$, the porosity $\varepsilon$, the mean particle radius $r$ and the specific interface area between electrolyte and active material $a$ \cite{Doyle1993}. Furthermore, \emph{effective} transport parameters need to be determined to account for the influence of the microstructure on the \emph{average} species transport.
While these models are able to describe the average battery behavior surprisingly well \cite{Tippmann2013,Doyle1996,Ecker2015b,Latz2015}, they cannot be expected to capture local microscopic effects. In particular, many degradation effects like, for instance, \Li plating depend on the local environment. Hence homogenized models cannot fully capture the interplay between microstructure and degradation phenomena with sufficient predictive power.
Therefore, more fundamental, spatially resolved models should be applied that are able to take the electrode microstructure explicitly into account \cite{Latz2011a,Latz2011b}. Without simplifications like volume averaging for the P2D-models these allow the computation of quantities on the scale of the electrode microstructure and are hence better suited for plating predictions.
To give an example, in \cite{Smith2009} a microstructure-based simulation study for a LCO-graphite battery was performed concentrating on the discharge behavior for a 2D cut of one given realization of the electrodes.
Although the numerical solution of these micro-models is computationally much more demanding they have been successfully applied to study cell performance \cite{Less2012a,Wang2007,Hutzenlaub2014}, coupling to thermal effects \cite{Latz2015,Yan2013,Taralov2014}, and to account for phase-separation dynamics within certain electrode materials \cite{Hofmann2016}. A framework for these spatially resolved simulations has been implemented in the software BEST \cite{ITWM2014}.

Lithium plating is one of the major degradation factors and security risks in \Liion batteries. Lithium plating describes the deposition of metallic \Li on the negative electrode \cite{Vetter2005}. This causes a loss of usable \Li (which reduces the cell's capacity) and might lead to the growth of \Li dendrites which can eventually create a short-circuit between the electrodes which can favor catastrophic thermal runaways. While model extensions to account for \Li plating are typically based on the porous electrode theory \cite{Tippmann2013,Newman1975,Arora1999,Tang2009,Legrand2014a} only very little work has been published where \Li plating models take the electrode microstructure into account \cite{monroe2003dendrite,akolkar2013mathematical}. In a recent publication \cite{Hein2016} a micro-scale model has been developed that is able to take the inhomogeneous electrode structure into account. 
The research presented in the current paper is based on this degradation model.

Spatially resolved electrochemical simulations as described above allow investigations of electrochemical behavior for realistic 3D microstructures. Thus, as input for these simulations, realistic 3D image data of battery electrodes is needed, which is already available even in-operando \cite{Finegan2016}. However, tomographic measurements of battery electrodes in 3D involve high costs and efforts. A methodology that has proven to be very promising in this context is stochastic microstructure modeling. Based on (only one or a few) tomographic measurements, a 3D parametric stochastic microstructure model can be constructed and calibrated using tools from stochastic geometry \cite{Chiu2013}. The model has been implemented in our software library GEOSTOCH \cite{Mayer2004}. It is able to describe the complex geometric microstructure in a statistical sense with only a few parameters such that each realization of the model represents the morphological characteristics of the tomographic image data (e.g., the distributions of particle size and shape, pore size distribution, etc.). Once fitted to tomographic image data, with hardly any effort an arbitrary number of realistic 3D microstructures can be generated on the computer. Moreover, systematic variation of model parameters allows the realization of virtual, but still realistic microstructures.
Such an approach has been considered, for example, in the context of organic solar cells \cite{SKTOJS_12}.
Using regression in the parameter space, microstructures that represent various manufacturing conditions could be generated on the computer and analyzed regarding their functionality. This results in an enormous reduction of complexity, as (most of) the structures do not have to be manufactured in the laboratory, but only tomographic image data of a few ones is needed. Similar examples of stochastic microstructure modeling can be found in literature \cite{Gaiselmann2013, Westhoff2015, Neumann2016}. In this work, a stochastic microstructure model for anodes in lithium-ion batteries \cite{Feinauer2015} is used.

While the aforementioned microscopic battery model can be solved by relatively standard iterative numerical methods, the solution process is computationally very demanding. 
In order to get meaningful results a sufficiently large electrode cutout needs to be resolved in the simulation. This results in huge time-dependent discrete systems which require considerable computing resources, already
for single simulation runs. Computational studies to identify critical parameters, to estimate the dependence 
of degradation on operating conditions or to support optimal design and control of batteries, however, require
many forward simulation runs with varying material or state parameters and are thus virtually impossible.
Hence, model reduction approaches for the resulting parameterized systems are indispensable for such simulation 
tasks. Concerning model reduction for \Liion battery models, we refer to the pioneering work \cite{CaiWhi:2009}
in the context of proper orthogonal decomposition (POD), and to the more recent contributions 
\cite{Cai2009,IlievLatzEtAl2012,VW13,LV15,OhlbergerRaveEtAl2014,OhlbergerRaveEtAl2016,OR16a} in the context of reduced basis methods.
In the work presented here, we rely on an implementation of recent model reduction methods
(such as the reduced basis method, POD, and the empirical interpolation method) 
implemented in our model order reduction library \PYMOR \cite{MRS16,PYMOR}.

The ability to efficiently and realistically predict the degradation behavior (here: \Li plating) of \Liion batteries under arbitrary load conditions relies on the following prerequisites:
\begin{enumerate}
  \item A physics-based predictive microscopic battery model that includes the plating mechanism.
  \item A method to create a number of virtual, yet realistic microstructures as basis to understand the correlation between structural properties and battery performance and degradation behavior.
  \item A numerical method that is able to efficiently perform a considerable number of three-dimensional, microstructure-resolving simulations for a variety of operating conditions.
  \item A software interface that is able to integrate these aspects into a common workflow.
\end{enumerate}
Within the project \MULTIBAT \cite{multibat} the authors developed and technically implemented a workflow that covers all the aforementioned aspects, namely stochastic geometry generation, model extension to account for plating, numerical implementation and development of model order reduction techniques. There are numerous papers in the literature (including several ones written by the authors of the present paper) on different components of the presented workflow. However, we are not aware of any publication on an algorithm integrating all these components into one single, holistic workflow, which enables comprehensive solutions of really complex problems related to Li-ion batteries. Thus, development, implementation, and testing of a holistic algorithm / workflow which integrates all components of the above-mentioned chain is one of the main contributions of the present paper. A feasibility study for a really complex problem, such as plating, is presented in order to illustrate the capabilities of the workflow that has been developed. The investigation of the interplay between 3D microstructure and electrochemical processes during the plating processes, which up to our knowledge has not been done so far in the literature, is another main contribution of the present paper.

The authors' developments on the individual components of the workflow have been reported earlier, details can be found in the references listed in the present paper, hence, these components are presented here relatively shortly. The emphasis in the present paper is given to the developed interfaces, to the integration of all the components into one single workflow, to the peculiarities related to the selected feasibility study, and to the parametric study of the plating process in stochastic geometry generation. Special attention is paid to computational efficiency, adapting the model reduction approach to the heavily nonlinear system of partial differential equations. The presented study reveals, that the complex information produced by the interplay between microstructure, lithium-ion transport and intercalation kinetics is hidden in a vastly reduced subspace of the full 3D information contained in the time-dependent scalar fields for lithium-ion concentration and electrochemical potential. The essential dynamics leading to plating in a complex microstructure can be represented by a sophisticated reduced model without losing spatial precision. The advantage of the model presented in this paper is the ability to perform fully 3D microstructure-resolved simulations of plating with nearly the same numerical efficiency as simulation with a P2D model i.e. a 1D volume averaged battery model, in which all structural details are lost.
We report on the \MULTIBAT workflow and briefly describe the details of all the individual aspects in Section \ref{sec:workflow}. In Section \ref{sec:results} we demonstrate the application of the developed methods by showing and discussing results of a simulation study and conclude with a summary in Section \ref{sec:summary}.

\section{The \MULTIBAT workflow}
\label{sec:workflow}
\begin{figure}
\centering
    \includegraphics{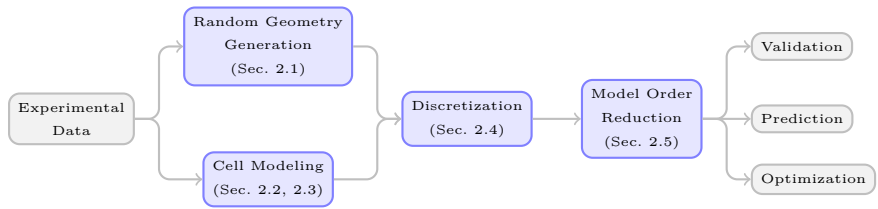}
	\caption{Schematic overview of the \MULTIBAT workflow.}\label{fig:abstract_workflow}
\end{figure}

In this section we discuss the individual components of the \MULTIBAT workflow (see Fig.~\ref{fig:abstract_workflow})
and their realization in more detail.
Based on experimental data, random electrode geometries with the same or modified structural
characteristics are generated (Sec.~\ref{sec:workflow:generation}), and a mathematical
model of the relevant electrochemical effects is formulated (Sec.~\ref{sec:workflow:cell_modeling} and
\ref{sec:workflow:degradation_modeling}).
The resulting continuum model is then discretized (Sec.~\ref{sec:workflow:discretization})
and reduced (Sec.~\ref{sec:workflow:MOR}), leading to a quickly computable microscale model of the
cell dynamics on realistic electrode geometries.
The software implementation and integration into a unified modeling and simulation workflow is
discussed in Sec.~\ref{sec:workflow:integration}.

It should be noted that, while we present a specific realization of the \MULTIBAT workflow targeted
at \Li plating, the same workflow can be applied to other questions in electrochemistry
and similar problem domains. Each individual component can be further developed and optimized
for other specific applications, independently of the other workflow components.

\subsection{Generation of random structures}
\label{sec:workflow:generation}
The study of local effects in the complex microstructures of battery anodes by spatially resolved models is computationally very expensive, particularly regarding random access memory. Therefore, only quite small sample sizes can be considered. As we are interested in local phenomena, there is a need for high-resolution of the images, which on the other hand means that the images typically only represent small cutouts of the material. This is why, in order to get reliable results, the computer experiments have to be carried out repeatedly using different samples. Furthermore, the imaging techniques are complex in preparation and involve long imaging times as well as high costs. This is why a suitable approach is to use randomly generated images of microstructures using parametric stochastic 3D models. This approach has already been used successfully in related applications for energy materials in fuel cells \cite{Gaiselmann2012} and solar cells \cite{SKTOJS_12}. A parametric stochastic model which describes the spatial structure is developed for the material and its parameters are fitted to image data. Using the calibrated model, an arbitrary number of structures that are similar to the image data in a statistical sense can be generated with hardly any effort. `Similar in a statistical sense' means that the realizations of the model do not resemble the image data exactly, but with respect to aggregated quantities and spatial properties. For example, simple characteristics like volume fraction and specific surface area can be matched, but also more complex spatial characteristics like the distribution of pore sizes or local tortuosity. Thus, realizations of a parametric stochastic microstructure model are an ideal input for spatially resolved electrochemical simulations. A further advantage is that their parameters can be changed to create virtual structures that have not been produced in the laboratory yet, and the electrochemical performance of those virtual structures can be analyzed on the computer, a procedure called virtual materials testing.

Here, we make use of a parametric stochastic 3D microstructure model for anode structures from \Liion battery cells \cite{Feinauer2015}. Besides the validation based on structural characteristics \cite{Feinauer2015} a validation using spatially resolved electrochemical simulations has been performed \cite{Hein2016b}. The variability of the modeling approach used here is demonstrated since the same model with some adaptions can be used to generate microstructures for energy cells \cite{Feinauer2015} and power cells \cite{Westhoff2016}.

\begin{figure}[!ht]
	\centering
	\subfigure[]{
		\includegraphics[width=0.22\textwidth]{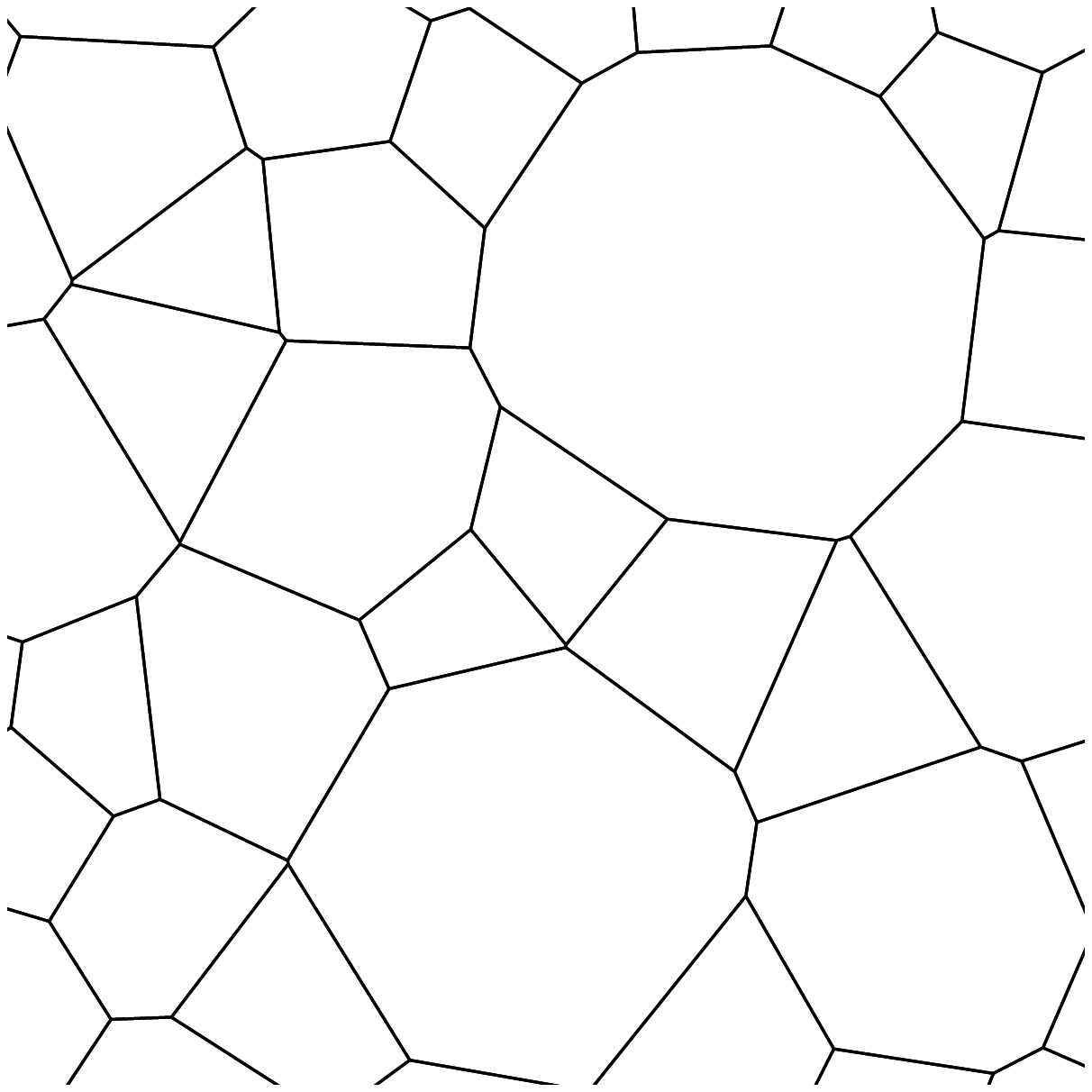}
		\label{fig:tessellation}
	}
	\subfigure[]{
		\includegraphics[width=0.22\textwidth]{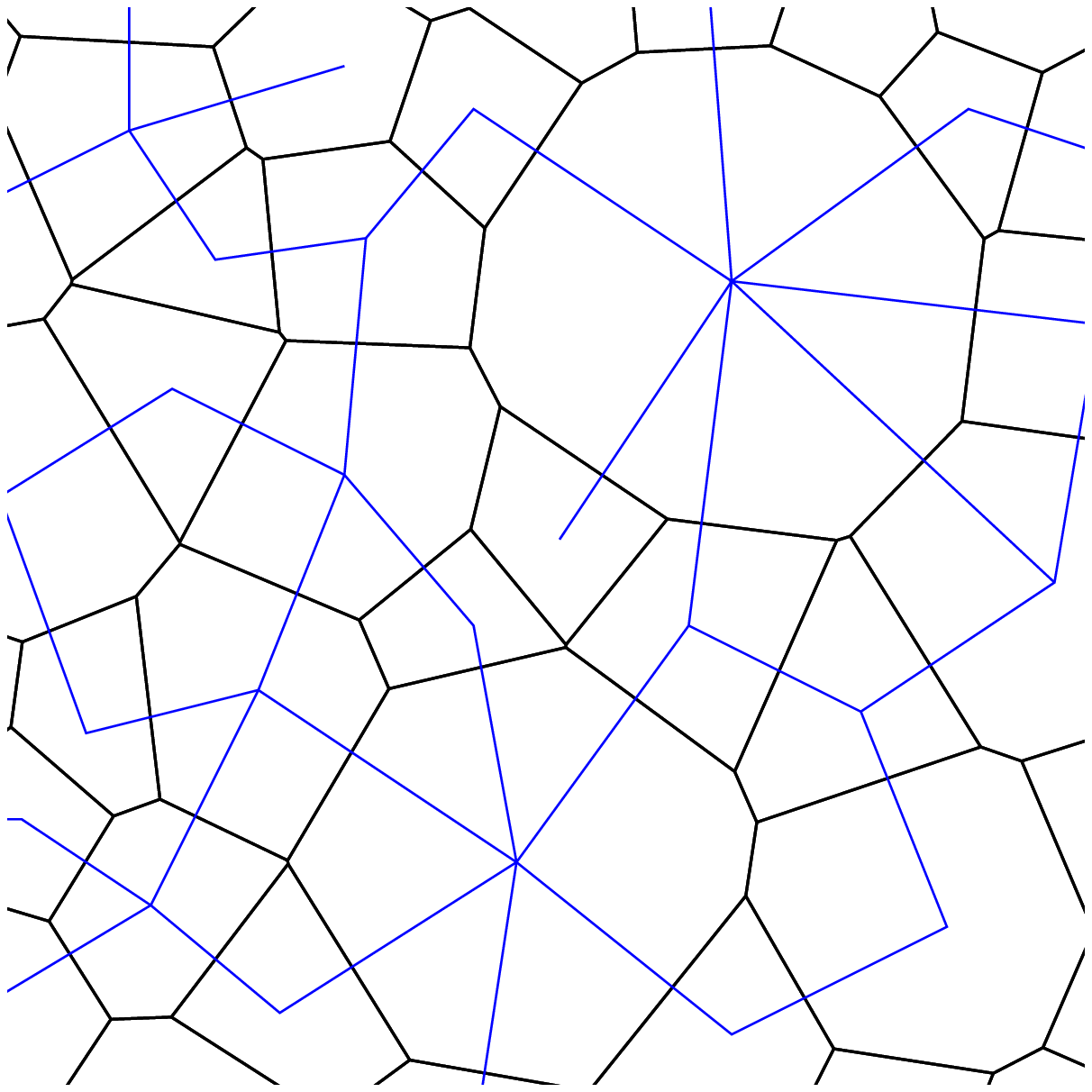}
		\label{fig:connectivity_graph}
	}
	\subfigure[]{
		\includegraphics[width=0.22\textwidth]{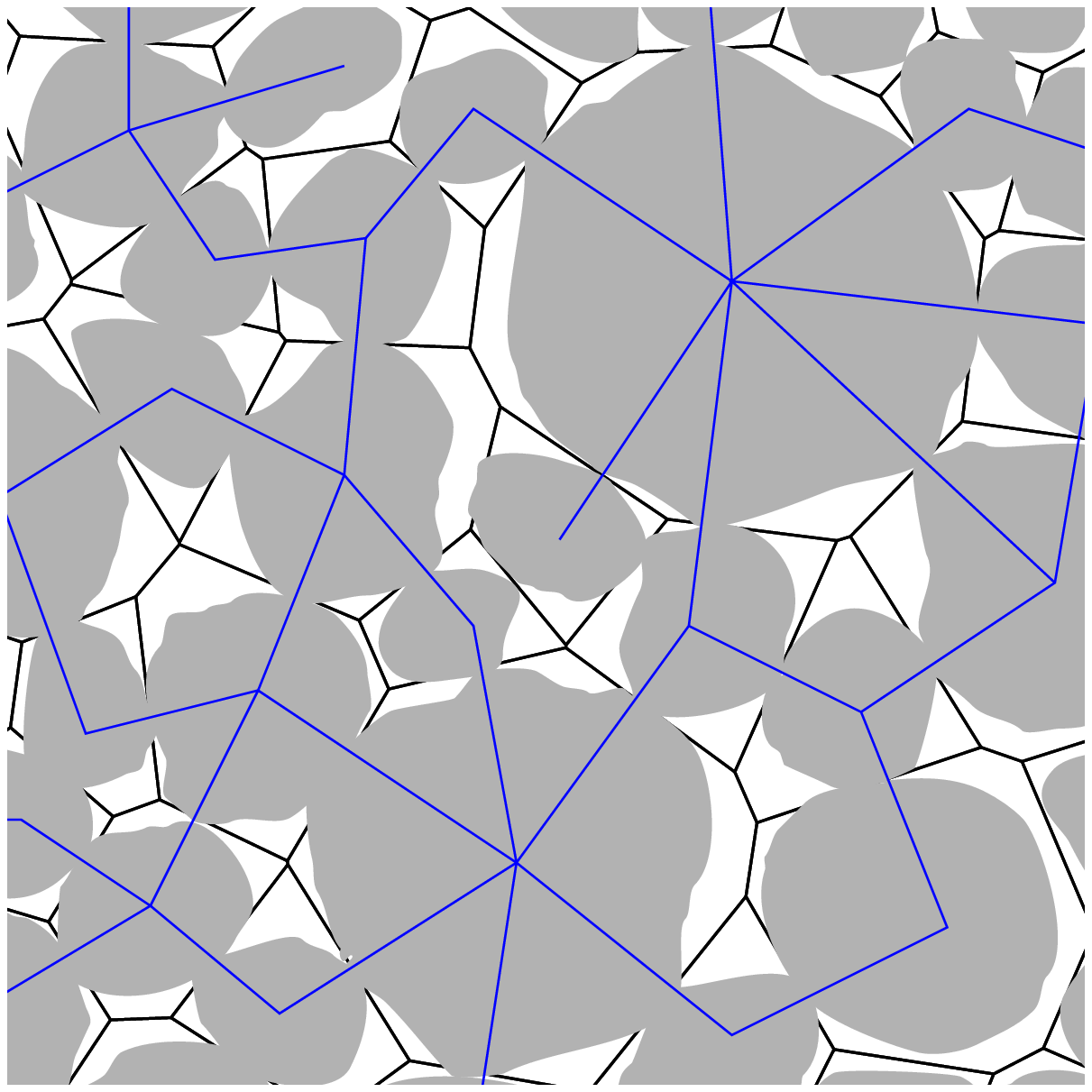}
		\label{fig:graph_and_particles}
	}\\
	\vspace{.1cm}
	\subfigure[]{
		\includegraphics[width=0.22\textwidth]{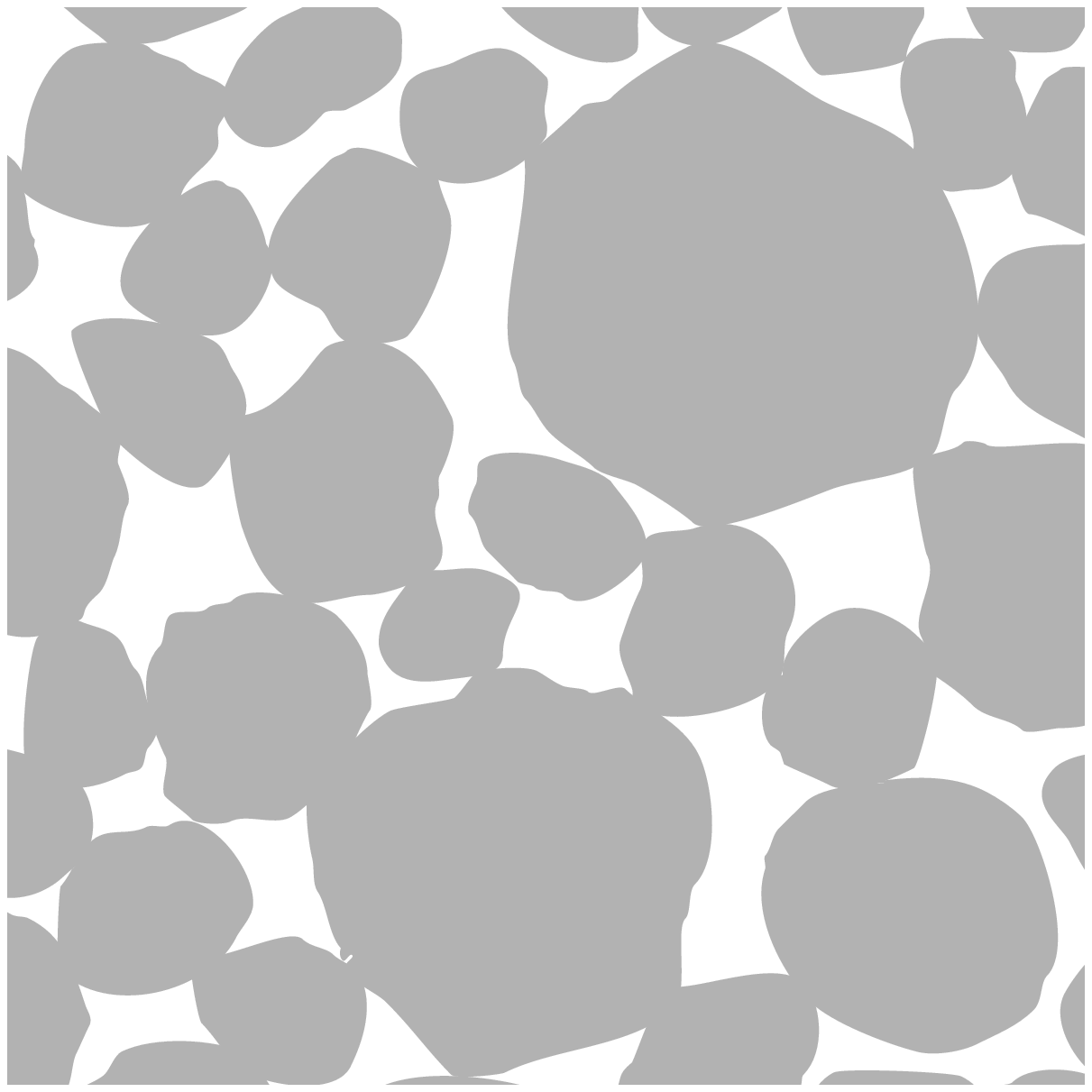}
		\label{fig:particles}
	}
	\subfigure[]{
		\includegraphics[width=0.22\textwidth]{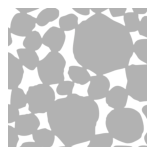}
		\label{fig:smoothed}
	}
	\caption{Schematic depiction of the stochastic model. (a) A random tessellation is produced, which roughly determines the particle shapes, sizes and locations. (b) A random graph describes how the particles are connected to each other. (c) The connected particles are generated using  random fields on the sphere. (d) and (e) The connected particles are retained and morphological smoothing is carried out. Reprinted from \cite{Feinauer2015} with permission from Elsevier.}
	\label{fig:stochastic_model}
\end{figure}

We now briefly recall some details of the stochastic 3D model that is used to generate the virtual anode microstructures used in the \MULTIBAT workflow. As mentioned above the model has already been published \cite{Feinauer2015} and all parameters as well as further details can be found there. Generally, the construction of the model consists of four steps that are also depicted in \figref{fig:stochastic_model}.

First, the locations, sizes and shapes of the particles are determined. Technically speaking, a Laguerre tessellation is generated (see~\figref{fig:tessellation}) based on a random sequential adsorption process. This tessellation decomposes the region of interest into convex polytopes. Later on, a particle is placed inside each of these polytopes. Thus, the Laguerre tessellation roughly indicates the spatial arrangement of particles. For details regarding tessellations, a broad spectrum of literature is available \cite{Chiu2013,Lautensack2007,Moller1989,Torquato2002}.

In the next step, a connectivity graph is generated that describes which particles are supposed to be connected, i.e., for each polytope $P$, we determine a set of neighboring polytopes $\{P_i, \,i=1,...,N\}, \, N\in\mathbb{N}$. The particles that are placed inside $\{P_i, \,i=1,...,N\}$ have to touch the particle in $P$. Full connectivity of all particles is ensured by the usage of a minimum spanning tree \cite{Prim1957}. Further connections are added to the minimum spanning tree depending on the size of the facet between two polytopes, as the probability of two particles being connected is larger for larger facet areas. Such a graph is depicted in \figref{fig:connectivity_graph}.

Now, a particle can be realized in each polytope, fulfilling the boundary conditions, i.e., touching the particles indicated by the connectivity graph. In more detail, the particles are modeled using Gaussian random fields on the sphere. Thus, the shape of the particles can be characterized by a mean radius $\mu$ and the angular power spectrum $A: [0, \infty) \rightarrow [0, \infty)$, see \cite{Lang2014}. The angular power spectrum is approximated by the function $A(l) = \frac{a l + b}{l^2 + cl + d}$ with coefficients $a = 0.4241$, $b = 0.356$, $c = -3.858$ and $d = 3.903$. In more detail, we do not use the mean radius $\mu$ directly but we generate the particles in a way that their volume is proportional to the volume of corresponding Laguerre cells.
The particles are sampled with the boundary conditions indicated by the connectivity graph using a special sampling algorithm that creates only realizations of the given Gaussian random field that fulfill those conditions.

The schematic depiction in \figref{fig:graph_and_particles} shows the particles with the tessellation and the connectivity graph. One can clearly see that the particles touch each other where indicated by the graph and on the other hand also fill their respective Laguerre polytopes. \figref{fig:particles} shows the system of connected particles without the tessellation and the connectivity graph as these are auxiliary tools that are no longer needed after the creation of particles.

Finally, a morphological smoothing \cite{Dougherty1992} is performed on the system of connected particles to mimic the effect of binder. In the given sample the volume fraction of the binder as well as the contrast in the tomographic images were too low to identify and model the binder as separate phase. From the known production process (slurry coating) we assume that this approach produces a similar effect as depicted in \figref{fig:smoothed}.

\begin{figure}[t]
	\centering
	\subfigure[]{
		\includegraphics[width=0.39\textwidth]{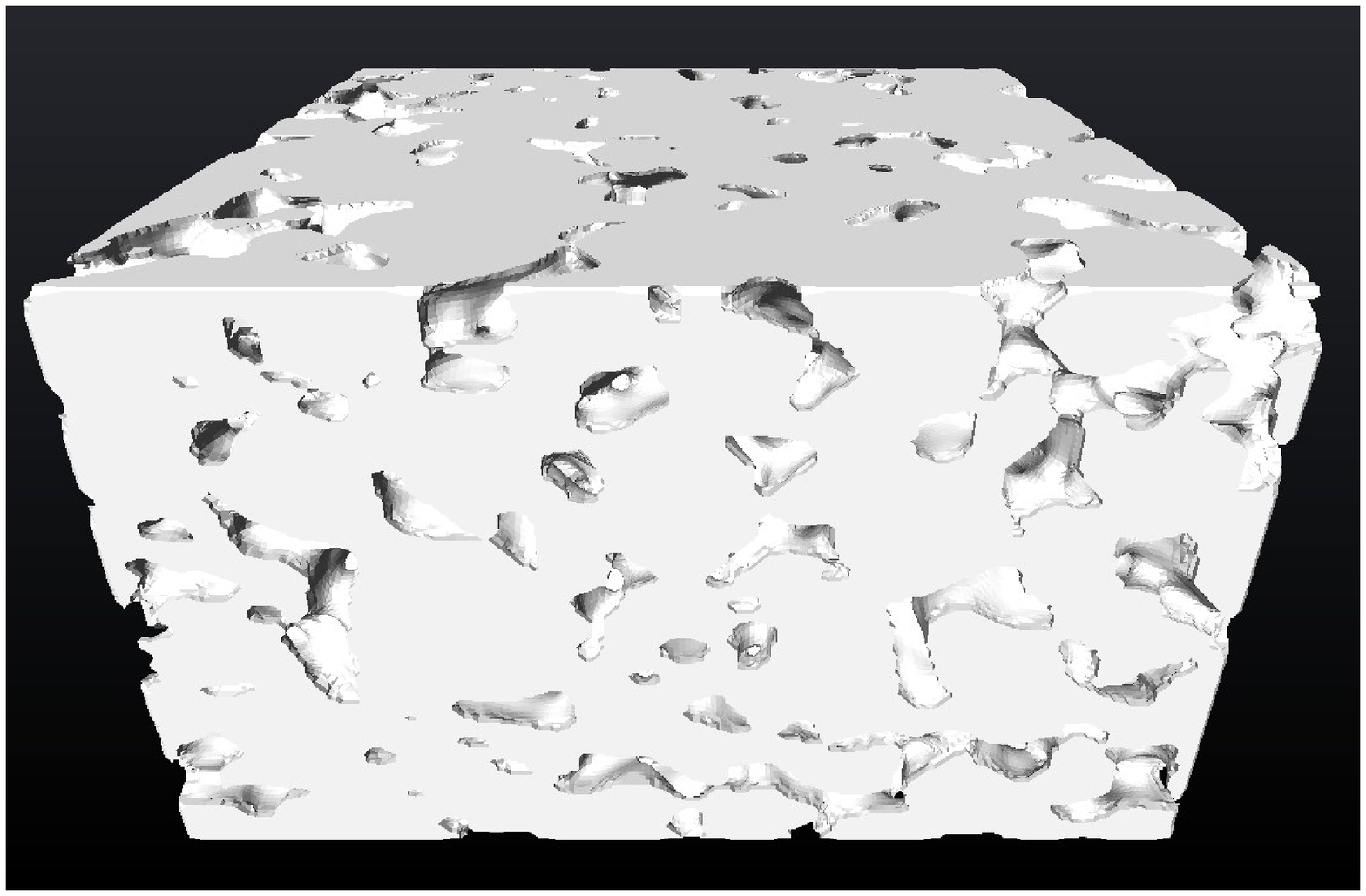}
		\label{fig:visual_comparison(a)}
	}
	\subfigure[]{
		\includegraphics[width=0.39\textwidth]{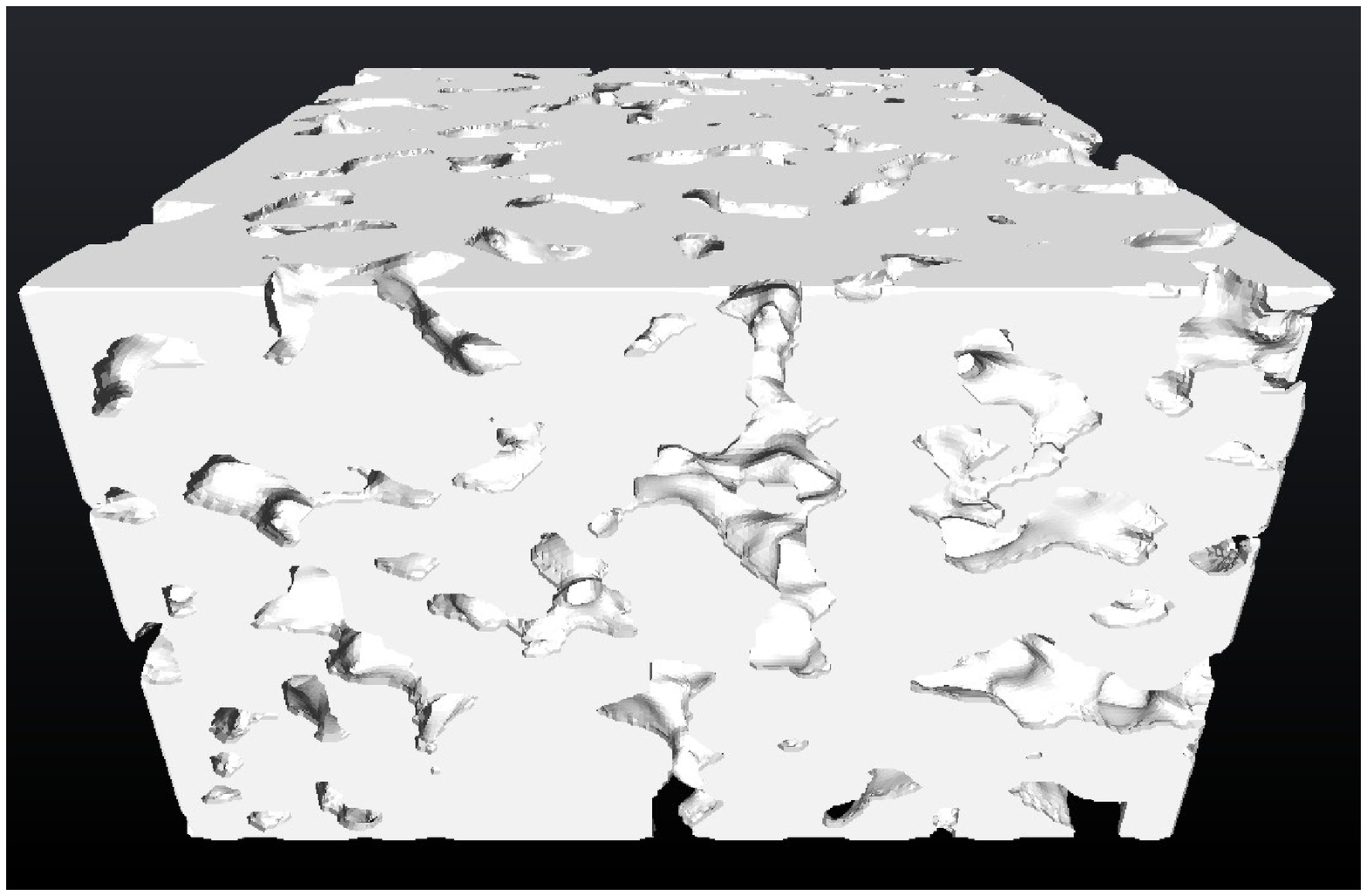}
		\label{fig:visual_comparison(b)}
	}
	\subfigure[]{
		\includegraphics[width=0.78\textwidth]{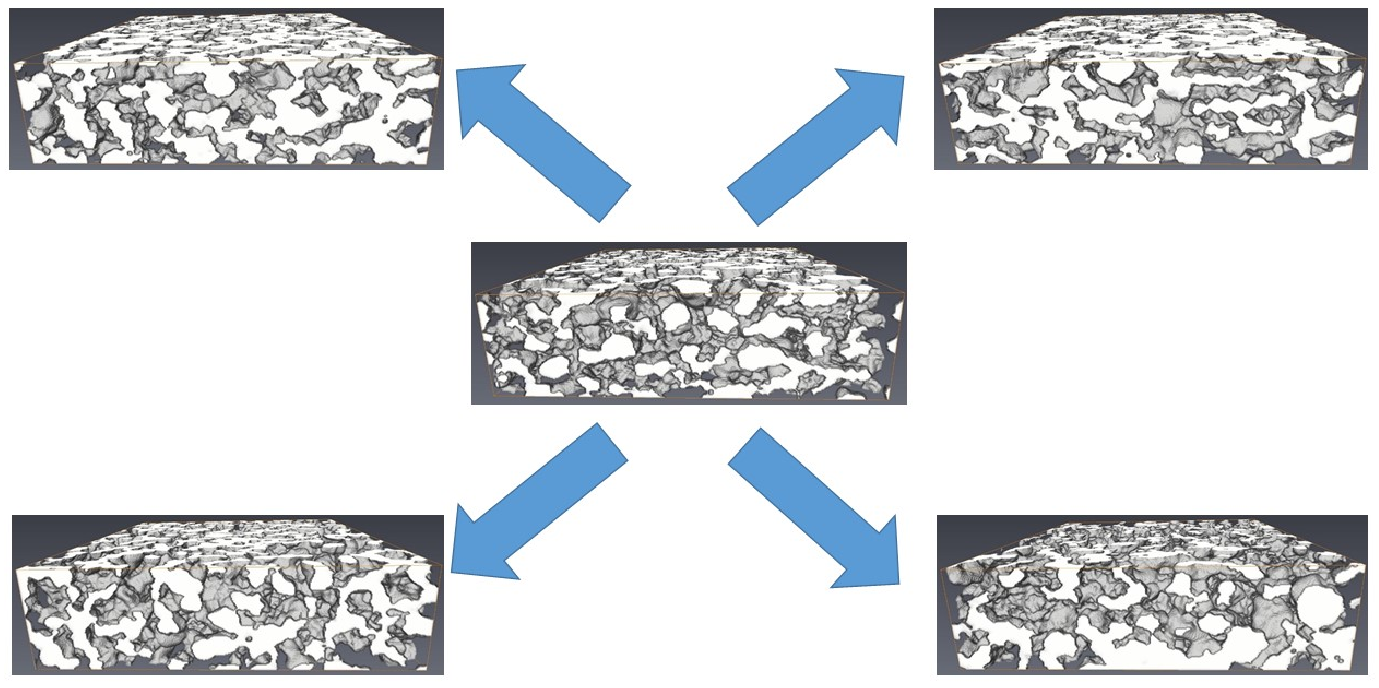}
		\label{fig:visual_comparison(c)}
	}	
	\caption{(a) 3D rendering of a cutout of a tomographic image of the energy cell anode. (b) 3D rendering of a simulated energy cell anode structure. (c) Realizations of the power cell model with various morphological properties; center: realization of the calibrated power cell model; top left: virtual structure with higher volume fraction of the particle phase; top right: virtual structure with more pronounced anisotropy effects; bottom left: virtual structure with no anisotropy effects; bottom right: virtual structure with decreasing volume fraction of the particle phase from bottom to top. Reprinted from \cite{Feinauer2015}(a+b) and \cite{Westhoff2016}(c) with permission from Elsevier.}
\end{figure}

In Fig. \ref{fig:visual_comparison(a)} and \ref{fig:visual_comparison(b)}, a cutout from the tomographic image data can be compared to a simulated anode structure. A very good visual accordance can be observed.

The model described so far is an excellent tool to generate virtual anode microstructures of energy cells, which are characterized by a high volume fraction of the solid phase. However, note that it can not directly be used to model the morphology of anodes in power cells, because due to the lower volume fraction of the solid phase, the boundary conditions for particles cannot be fulfilled reasonably. Therefore, an extension of the model has been proposed \cite{Westhoff2016}. To account for the lower volume fraction, a Laguerre tessellation with marked polytopes is used. The polytopes are marked either as `filled', i.e., a particle is placed in the polytope, or as `empty', which means that no particle is placed here. Thereby, a reasonable allocation of the different polytopes as well as full connectivity of the resulting structure is ensured. Furthermore, the model is able to include anisotropy effects of the solid phase, i.e., particles can be elongated in horizontal direction rather than in vertical direction. This results in a remarkable flexibility such that the model can be used to create a broad spectrum of virtual anode microstructures with a variety of morphological properties, see Fig. \ref{fig:visual_comparison(c)} for some examples.

In this study, we focus on electrochemical simulations of anode microstructures in energy cells. Thus, all microstructures which are discussed in the present paper, are created using the energy cell model \cite{Feinauer2015}. Electrochemical simulations on virtual structures generated by the power cell model are subject of further research.

\subsection{Microscopic cell modeling}
\label{sec:workflow:cell_modeling}
In this section we recollect the equations that describe the transport of \Liions in a three-dimensional microstructure generated by the method as described in the previous section.
The physical model has been derived based on species, charge and energy conservation to yield a set of equations that describe the spatial and temporal distribution of \Liions, electrical potentials and temperature \cite{Latz2011b,Latz2013}.
However, in order to fit to the isothermal plating model we neglect the effects of heat production here and consider an isothermal system where temperature enters as a constant model parameter.

For the purpose of this study we restrict ourselves to half-cell simulations, i.e. we consider a porous graphite electrode modeled by the discussed stochastic method against a \Li foil as counter-electrode. The simulation domains are connected with the external operation conditions through dedicated current collector phases on each electrode.
The remaining space of the computational domain is filled with ion conductive electrolyte.

Within the graphite particles we have the following equations for \Li concentration $\var{c}{Gr}$ and electrical potential $\var{\Phi}{Gr}$

\begin{alignat}{2}
    \partial_t \var{c}{Gr} &= -\bm{\nabla} \cdot \bm{\var{N}{Gr}}
                           &&=-\bm{\nabla} \cdot \left[-\var{D}{Gr} \bm{\nabla}  \var{c}{Gr} \right] \label{eq:cgr}\,, \\
    0 & = -\bm{\nabla} \cdot \bm{\var{j}{Gr}}
      &&= -\bm{\nabla} \cdot \left[ -\var{\sigma}{Gr} \bm{\nabla} \var{\Phi}{Gr}\right]\,, \label{eq:phigr}
\end{alignat}
where $\var{D}{Gr}$ is the \Li diffusion coefficient and $\var{\sigma}{Gr}$ is the electrical conductivity of the material. The ion flux and electric current density are denoted by $\bm{\var{N}{Gr}}$ and $\bm{\var{j}{Gr}}$, respectively.
Also in the domains of the \Li foil and the current collectors electronic conduction is considered and hence \eqref{eq:phigr} is also solved in these domains (with the respective conductivities of course).
Since there is no intercalation and diffusion of ions neither in the \Li counter-electrode nor in the current collectors, \eqref{eq:cgr} is only relevant for the graphite domain.

Within the electrolyte domain, \Li concentration $\var{c}{El}$ and electrochemical potential $\var{\varphi}{El}$ are coupled through

\begin{alignat}{2}
\partial_t \var{c}{El} &= -\bm{\nabla} \cdot \bm{\var{N}{El}}
		       &&= -\bm{\nabla} \cdot \left[ -\var{D}{El} \bm{\nabla}  \var{c}{El} + \frac{\var{t}{+}}{F} \var{\bm{j}}{El}\right]\label{eq:cel}\,,\\
0 & = -\bm{\nabla} \cdot \bm{\var{j}{El}}
  &&= -\bm{\nabla}\cdot \left[- \var{\kappa}{El}\bm{\nabla} \var{\varphi}{El} - \var{\kappa}{El} \frac{\var{t}{+}-1}{F}  \frac{\partial\mu }{\partial\var{c}{El}}  \bm{\nabla} \var{c}{El}\right] \,, \label{eq:phiel}
\end{alignat}
where $t_+$ is the transference number of lithium in the electrolyte, $\var{\kappa}{El}$ is the ionic conductivity of lithium inside the electrolyte and $F$ is the Faraday constant. The derivative of the electrolyte chemical potential is given as

\begin{align}
\frac{\partial\mu }{\partial\var{c}{El}} &= \frac{R\cdot T}{\var{c}{El}} \cdot \left( 1 + \frac{\partial\log{\var{f}{+}}}{\partial\log{\var{c}{El}}}\right)\,,
\end{align}
with $T$ denoting the temperature, $R$ the gas constant and $f_+$ the activity coefficient.

On the interfaces between the electrodes and electrolyte two types of reactions need to be considered: That is an intercalation reaction on the graphite side and a \Li deposition reaction on the counter-electrode side. The different phases (graphite, electrolyte and counter-electrode) are coupled via interface conditions

\begin{equation}\label{eq:interfaceCond}
\begin{split}
\var{\bm{j}}{El}\cdot \var{\bm{n}}{So-El} &= \var{i}{interface}\,,\\
\var{\bm{j}}{So}\cdot \var{\bm{n}}{So-El} &= \var{i}{interface}\,,\\
\var{\bm{N}}{El}\cdot \var{\bm{n}}{So-El} &= \frac{\var{i}{interface}}{F}\,,\\
\var{\bm{N}}{So}\cdot \var{\bm{n}}{So-El} &= \frac{\var{i}{interface}}{F}\,,
\end{split}
\end{equation}
with ``So'' (solid) being either graphite or metallic \Li. By convention the interface normal \var{\bm{n}}{So-El} points from solid into the electrolyte.
These conditions express the continuity of the current and mass fluxes between the phases. The current flow through these interfaces depends on the corresponding reactions. The intercalation reaction is described by a Butler-Volmer-like expression \cite{Latz2013}

\begin{align}
 \var{i}{interface} &= \var{i}{Gr-El} \notag\\&= 2\cdot \vart{i}{Gr-El}{00} \cdot \sqrt{\var{c}{Gr} \cdot \var{c}{El}} \cdot \sinh \left(\frac{F}{2\cdot R \cdot T}  \cdot \var{\eta}{Gr-El}\right)\,,\label{eq:ise_gr}
\end{align}
where the overpotential is given by $\var{\eta}{Gr-El} = \var{\Phi}{gr} - \vart{U}{0}{Gr} - \var{\varphi}{El}$. The electrode's open-circuit potential $\vart{U}{0}{Gr}$ is a concentration dependent material property. The rate constant $\vart{i}{Gr-El}{00}$ depends on the lithium salt and the electrolyte composition. The transfer coefficients $\alpha_c+\alpha_a =1$ of the intercalation reaction were assumed to be symmetrical ($\alpha_{a,c}=0.5$). The form of \eqref{eq:ise_gr} differs from the usual Butler-Volmer expression by omitting the common $\left(\vart{c}{Gr}{max}-\var{c}{Gr}\right)^\alpha_a$ term, since a rigorous thermodynamically consistent derivation does not in general yield this prefactor \cite{Latz2013}. 
The Butler-Volmer model and the used exchange current include a relation between the potentials and the current flux. The exponential shape of the exchange current introduces highly pronounced non-linearities into the numerical system.
The reaction at the counter-electrode is described by a simple exchange current, which minimizes the effect of the counter-electrode on the simulation results

\begin{align}
    \var{i}{interface}= \var{i}{CE-El} &= 2\cdot \vart{i}{CE-El}{00} \cdot \sinh \left(\frac{F}{2\cdot R \cdot T}\cdot \var{\eta}{CE-El} \right)\,,\label{eq:ise_li}
\end{align}
where the rate constant is given by $\vart{i}{CE-El}{00}$.
On the \Li electrode the overpotential is simply given by $\var{\eta}{CE-El}=\var{\Phi}{CE} - \var{\varphi}{El}$.
The remaining interface conditions are shown in Tab.~\ref{tab:interfaceCond}.

\begin{table*}
    \caption{Overview of interface conditions between the different material domains for the ion fluxes $N$ and current densities $j$. ``cont.'' mathematically means no interface but a continuous flux according to the transport equation. Between the graphite electrode and the \Li foil as counter electrode there is obviously no interface.\label{tab:interfaceCond}}
{\small
\medskip
\scriptsize
\begin{tabular} {r|ccccc}
              & graphite & electrolyte & \Li foil & current collector & plated \Li\\\hline
              graphite      & cont.    & \eqref{eq:interfaceCond},\eqref{eq:ise_gr} & no interface & $N=0, j=\text{cont.}$ & $N=0,j=\text{cont.}$\\
              electrolyte   & \eqref{eq:interfaceCond},\eqref{eq:ise_gr} & cont. & \eqref{eq:interfaceCond},\eqref{eq:ise_li} & $N=-1, j=0$ & \eqref{eq:interfaceCond},\eqref{eq:plelinterface}\\
              \Li foil      & no interface & \eqref{eq:interfaceCond},\eqref{eq:ise_li} & cont. & $N=0, j=\text{cont.}$ & no interface\\\
\!\! current col.  & $N=0, j=\text{cont.}$\!\!\!   & $N=0, j=0$\!\!\!   & $N=0, j=\text{cont.}$ \!\!\!\!\! & cont. &  $N=0, j=\text{cont.}$\!\!\!\\\
\!\! plated \Li  & $N=0,j=\text{cont.}$   & \eqref{eq:interfaceCond},\eqref{eq:plelinterface}  & no interface & $N=0, j=\text{cont.}$\!\!\! & cont.
\end{tabular}
}
\end{table*}

These equations describe an \emph{ideal} battery, i.e. no degradation processes are considered. This extension is outlined in the next section.

\subsection{Electrochemical degradation modeling}
\label{sec:workflow:degradation_modeling}
The focus in this work is the degradation process \textit{\Li plating}, where the \Liions form an unwanted metallic phase on the surface of the intercalation material of the negative electrode. The electrochemical modeling of this process is briefly described in the following section, details can be found in the corresponding publication \cite{Hein2016}. Two states of the \Liions are of relevance: the \Liions dissolved in the electrolyte $\chem{Li}{El}{+}$ and the metallic/plated \Li phases $\chem{Li}{Pl}{0}$. The transition between these two is expressed by the reaction

\begin{align}
\chem{Li}{Electrolyte}{+} + \chem{e}{Electrode}{-} \rightleftharpoons \chem{Li}{Plated}{0}\,. 
\label{eq:ap2:theo:reaction:plating}
\end{align}
The overpotential of the plating and stripping reaction is defined by the difference between the electrochemical potential $\tilde{\mu}$ of the two \Li phases involved\cite{Newman2004}

\begin{align}
F \cdot \var{\eta}{Pl/St} = \vart{\tilde{\mu}}{Li^+}{Pl} - \vart{\tilde{\mu}}{Li^+}{El}\,. \label{eq:ap2:theo:overpot:eq1}
\end{align}
With the definition of the reference state $\vart{\mu}{Li^0}{Pl} = \vart{\mu}{Li^+}{Pl} + \vart{\mu}{e^-}{Pl}$ and the electrochemical potentials of \Liions inside the electrolyte and a solid phase (see \cite{Hein2016}), the overpotential \eqref{eq:ap2:theo:overpot:eq1} can be rewritten to

\begin{align}
\var{\eta}{Pl/St} = \var{\Phi}{Pl} - \vart{\varphi}{Li^+}{El}\,, \label{eq:ap2:theo:overpot:eq2}
\end{align}
where $ \vart{\varphi}{i}{p}$ denotes the electrochemical potential of species $i$ in phase $p$ with respect to the reference state $\vart{\mu}{Li^0}{Pl}$.
Plating of \Li is occurring if the overpotential reaches negative values ($\var{\eta}{Pl/St} < 0$). The metallic \Li phase on the surface of the anode is not in a stable configuration, even if no external current is applied to the system. As soon as \Li is plated on the surface of the active material, the \Li metal can react with its surroundings in different ways. The reaction between the plated \Li and the electrolyte results in the growth of a solid-electrolyte interphase (SEI), which leads to an irreversible loss of lithium \cite{Vetter2005}. Apart from phenomenological models no theory exists which combines lithium intercalation, lithium plating and SEI growth. Hence, this irreversible pathway is not included in the present paper. The plated \Li can also intercalate charge-neutrally into the supporting graphite. This reaction represents a reversible \Li stripping pathway. We are not aware of any literature regarding the identification and parameterization of the charge-neutral reintercalation. Hence the direct reintercalation from the plated lithium into graphite was neglected in this work. 

The stripping and plating reaction of the \Li is described by a Butler-Volmer-like equation

\begin{align}
  \var{i}{Pl-El} = &\vart{i}{Pl-El}{00} \cdot \sqrt{\var{c}{El}}\notag\\ &\cdot \left(\var{f}{pre}\left(\var{n}{Li}\right) \cdot \exp\left(\frac{ F\var{\eta}{Pl-El}}{2\cdot R\cdot T} \right) - \exp\left(-\frac{F\var{\eta}{Pl-El} }{2\cdot R\cdot T} \right) \right)\,.\label{eq:iplel}
\end{align}
The Bulter-Volmer-like expression is derived for non-vanishing phases. But, the plated \Li phase can completely desolve during stripping. Hence, the vanishing of the plated \Li phase is considered in the exchange current by the numerical regularization function $\var{f}{pre} \left( \var{n}{Li} \right)$, which depends on the amount of plated \Li $\var{n}{Li}$

\begin{align}
\var{f}{pre}\left(\var{n}{Li}\right) &= \frac{\left(\var{n}{Li}\right)^4}{\left(\vart{n}{Li}{const}\right)^4 + \left(\var{n}{Li}\right)^4}\,.
\end{align}
Based on numerical considerations \cite{Hein2016}, we set the constant $\vart{n}{Li}{const}$ to a value corresponding to a thickness of plated lithium of $0.48\,\text{nm}$. For partially covered surfaces more detailed models are necessary to capture the stripping of partially covered surfaces including the surface-driven dissolution of small lithium droplets.

At the interface between the plated \Li and the electrolyte the current through the interface is equal to the stripping current

\begin{align}
  \var{i}{interface} &= \var{i}{Pl-El}\,.\label{eq:plelinterface}
\end{align}
All the interface conditions which are relevant for electrochemical simulations in this paper are listed in Tab.~\ref{tab:interfaceCond}.

\begin{figure}[!ht]
	\centering
	\subfigure[]{
		\includegraphics[width=0.4\textwidth]{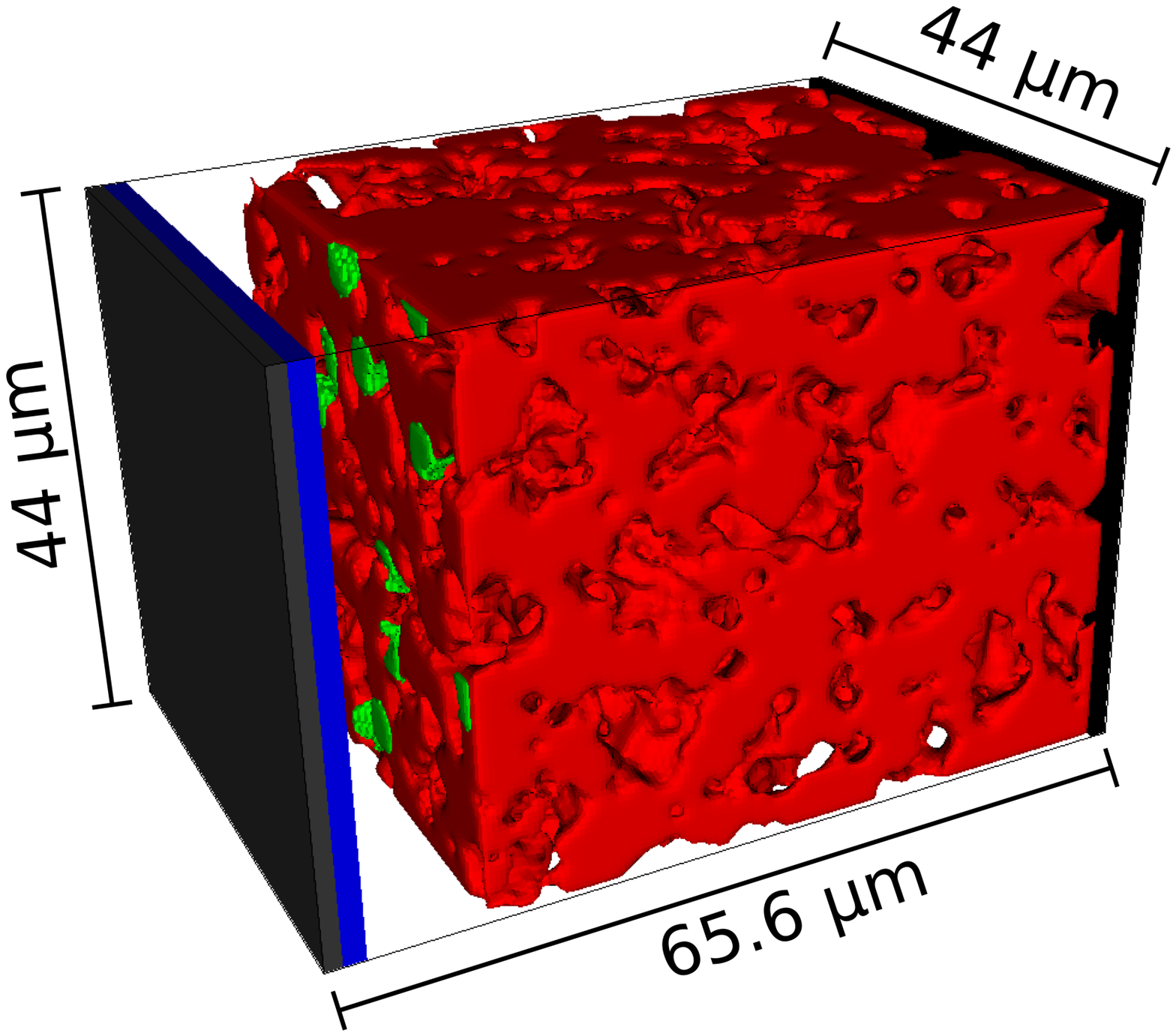}
		\label{fig:microsrtucture}
	}
	\subfigure[]{
		\includegraphics[width=0.4\textwidth]{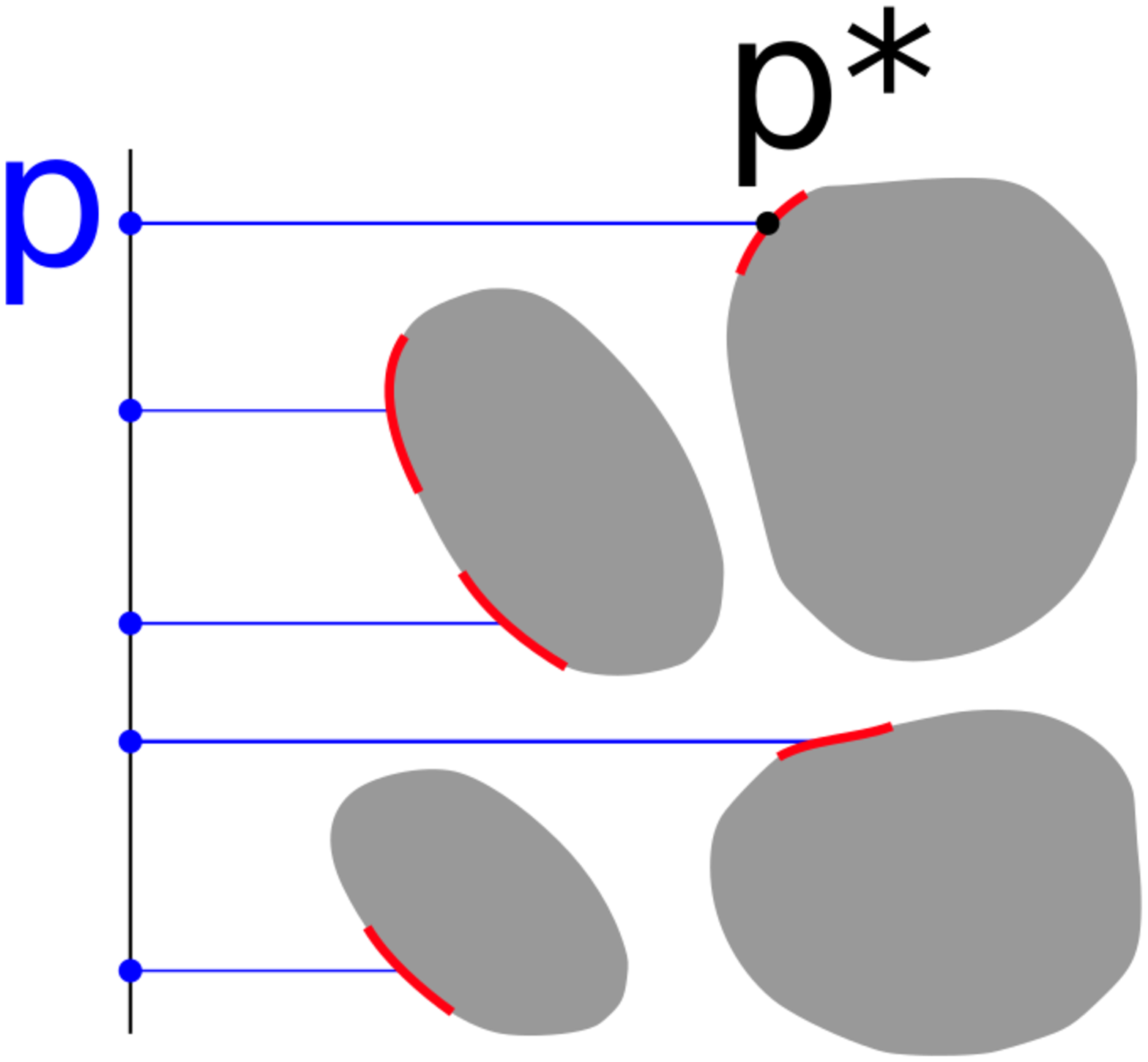}
		\label{fig:plating_model}
	}
 \caption{(a) Example of a 3D microstructure generated based on the stochastic simulation algorithm as described in \secref{sec:workflow:generation}. The plated \Li is shown as green spots on the separator-graphite interface. (b) Schematic depiction of the simulated initial \Li plating in 2D. The blue dots indicate the starting positions of the grains. The red lines indicate the simulated \Li plating.}
\end{figure}

In this paper the stripping process of plated \Li is simulated by including the plated \Li into the 3D microstructure as an additional volume phase. In \figref{fig:microsrtucture} an example of a 3D microstructure with plated \Li is shown. The porous electrode (red/right) is generated by the stochastic generation algorithm as described in \secref{sec:workflow:generation}. Additionally, regions with plated \Li are positioned randomly at the separator-electrode interface.  
The microstructure shown in Fig. \ref{fig:microsrtucture} is used for the model order reduction experiment described in \secref{sec:results:mor}.

One important deviation from the microstructure model \cite{Feinauer2015, Westhoff2016} was to introduce a third phase in the anode structure, plated \Li. We use a fairly simple model to create a slightly plated structure as initial condition. These initial conditions are, like the microstructure model, simulated stochastically. This means that for every run the microstructure model extended by the plated \Li phase generates a new structure but with similar statistical properties.

The method used to create the plated \Li phase is a germ-grain model. This means that in a first step germs are simulated and in the second step grains are placed around the germs \cite{Chiu2013}. The parameters of this model are the intensity $\lambda$ of the Poisson process that is used to generate the germs and the grain radius $r_s$. Due to the lack of experimental data the values where chosen to $\lambda = 0.01/\si{\micro \metre}^2$ and $r_s = 2.2\,\si{\micro\metre} $.
The general idea of the model is to place plating germs randomly on the surface of the particles in the electrode and then initialize plating around the germs on all the points on the structure's surface within the radius given by $r_s$. In more detail the following is done, see Fig. \ref{fig:plating_model}:
\begin{itemize}
	\item Select points $\{p_i, i\in\mathbb{N}\}$ at the separator-anode interface via a Poisson point process with intensity $\lambda$.
	\item For each point $p_i, i\in\mathbb{N}$ do the following:
	\begin{itemize}
		\item Find the first point $p_i^*$ where the straight line from $p_i$ towards the anode current collector interface meets the anode structure. Let $\Theta$ be the particle on the boundary of which $p_i^*$ is placed.
		\item Consider a sphere $B(p_i^*, r_s)$ around $p_i^*$ with radius $r_s$. Let $L_i=\{x\in\mathbb{R}^3: x\in B(p_i^*, r_s) \text{ and } x\in \partial\Theta\}$, where $\partial\Theta$ is the boundary of $\Theta$.
	\end{itemize}
	\item Let $L=\bigcup\limits_{i\in\mathbb{N}} L_i$ be the plated lithium phase.
\end{itemize}
The plated lithium phase $L$ is then discretized as a one voxel thick phase on the surface of the particles. The material parameters and reaction constants of the ion transport and plating model are adapted from a previous publication \cite{Hein2016}.

\subsection{Discretization and high-dimensional simulation}
\label{sec:workflow:discretization}
For the spatial discretization of the presented plating model, a cell-centered finite volume scheme
on a uniform voxel grid is considered. Hence the discretization is naturally conservative. Furthermore,
with the simple grid structure meshing of the complex microstructure is straight-forward.
The Butler-Volmer interface conditions \eqref{eq:interfaceCond} are
prescribed as the numerical flux across the respective domain interfaces, leading to a
global space differential operator on the entire computational domain.
Choosing implicit Euler time stepping for time discretization, we obtain a
series of discrete nonlinear equation systems of the form

\begin{align}
\label{eq:ap4:detailed}
        \begin{bmatrix}
		\frac{1}{\Delta t^{(n+1)}}(c_{\MU}^{(n+1)} - c_{\MU}^{(n)}) \\
                 0
        \end{bmatrix}
         + A_\MU
         \left(\begin{bmatrix}
                c_{\MU}^{(n+1)} \\
                \varphi_{\MU}^{(n+1)}
         \end{bmatrix}\right)
         = 0, \notag\\ (c_{\MU}^{(n+1)}, \varphi_{\MU}^{(n+1)}) \in V_h \oplus V_h.
\end{align}
Here, $V_h$ denotes the discrete finite volume space of locally constant grid
functions, $c_\MU^{(n)}, \varphi_\MU^{(n)} \in V_h$
denote the concentration and potential fields at time step $n$ for some
$p$-tuple of model parameters $\MU$ contained in a parameter domain of interest $\mathcal{P} \subset \mathbb{R}^p$, and
$A_\MU: V_h \oplus V_h \to V_h \oplus V_h$ is the finite volume space differential operator.
The system is closed by $c^{(0)}_\MU = c_0$ for some fixed initial \Li distribution
$c_0 \in V_h$.
The time step size $\Delta t^{(n)}$ is chosen adaptively for each time step to
accommodate the different time scales during and after the stripping of the plated \Li.
The nonlinear equation systems are solved in \BEST using Newton's method and an
algebraic multigrid solver for the solution of the linear correction equations.
Details on the discretization are provided in a previous publication \cite{Popov2011}.

\subsection{Model order reduction and reduced simulation}
\label{sec:workflow:MOR}

The computation of a single solution trajectory $c_\MU^{(n)},\varphi_\MU^{(n)}$ requires many
hours, even for relatively small geometries (cf.~\secref{sec:results:mor}).
In order to make parameter studies computationally feasible, reduced basis model reduction
techniques\cite{QuarteroniManzoniEtAl2016,HesthavenRozzaEtAl2016,Ha14} are applied
which have been implemented in our model order reduction library \PYMOR \cite{MRS16,PYMOR}.
This allows us to obtain a quickly solvable reduced order model
approximating the full order model \eqref{eq:ap4:detailed}.

To construct the reduced order model, solutions of \eqref{eq:ap4:detailed} are computed for few appropriately selected
parameters $\MU_1, \ldots, \MU_S$.
Various advanced algorithms exist for the selection of these \emph{snapshot} parameters,
often based on a greedy search procedure (cf. \cite{QuarteroniManzoniEtAl2016,HesthavenRozzaEtAl2016,Ha14}).
In our basic test case with a one-dimensional parameter domain $\mathcal{P}$ (\secref{sec:results:mor}), a
simple equidistant parameter sampling will be sufficient, however.

From this data, reduced approximation spaces $\tilde{V}_c, \tilde{V}_\varphi$ for the concentration
and potential fields are constructed
via proper orthogonal decomposition (POD, principal component analysis) \cite{sirovich87} of the snapshot data
sets $\mathcal{S}_c = \{c^{(n)}_{\MU_s}, c^{(n, i)}_{\MU_s}\}$,
$\mathcal{S}_\varphi = \{\varphi^{(n)}_{\MU_s}, \varphi^{(n, i_k)}_{\MU_s}\}$.
Here, $c^{(n,i)}_{\MU_s}$, $\varphi^{(n,i)}_{\MU_s}$ denote the intermediate Newton stages
during the solution of \eqref{eq:ap4:reduced}, which are included for improved numerical stability.
By construction, we in particular have $\tilde{V}_c \subseteq \Span \mathcal{S}_c$ and
$\tilde{V}_\varphi \subseteq \Span \mathcal{S}_\varphi$.
While $\dim V_h$ is in the order of $10^6$, we typically have $\dim \tilde{V}_c, \dim \tilde{V}_\varphi < 100$,
which makes significant computational speedups possible.

After the reduced approximation space $\tilde{V} = \tilde{V}_c \oplus \tilde{V}_\varphi$ has
been computed, the reduced order model is obtained via Galerkin projection of \eqref{eq:ap4:detailed} onto
$\tilde{V}$. I.e., we solve

\begin{align}
\label{eq:ap4:reduced}
		P_{\tilde{V}} \left\{
                \begin{bmatrix}
			\frac{1}{\Delta t^{(n+1)}}(\tilde{c}_{\MU}^{(n+1)} - \tilde{c}_{\MU}^{(n)}) \\
                         0
                \end{bmatrix} 
                 + 
                 A_\MU
                 \left(\begin{bmatrix}
		        \tilde{c}_{\MU}^{(n+1)} \\
			\tilde{\varphi}_{\MU}^{(n+1)}
                 \end{bmatrix}\right) \right\}
         = 0, \notag \\ (\tilde{c}_{\MU}^{(n+1)}, \tilde{\varphi}_{\MU}^{(n+1)}) \in
	 \tilde{V},
\end{align}
$\tilde{c}^{(0)}_\MU=P_{\tilde{V}_c}(c_0)$, where $P_{\tilde{V}}$ / $P_{\tilde{V}_c}$
denotes the $L^2$-orthogonal projection onto $\tilde{V}$ / $\tilde{V}_c$.

However, even though \eqref{eq:ap4:reduced} contains only $\dim \tilde{V}$
degrees of freedom, its solution requires the evaluation of the high-dimensional
system operator $A_\MU$.
This strongly limits the achievable speedup in computation time
when solving \eqref{eq:ap4:reduced} instead of \eqref{eq:ap4:detailed}.

To overcome this issue, $A_\MU$ is replaced by a quickly evaluable low-order
approximation using the empirical interpolation technique \cite{BarraultMadayEtAl2004,DrohmannHaasdonkEtAl2012}:
for an arbitrary (nonlinear) operator $O: X \to Y$, the \textsc{EI-Greedy} algorithm
is used to compute a low-order interpolation
space $\tilde{Y} \subseteq Y$ from evaluations of $T$ on given solution
trajectories, after which the interpolated operator $\mathcal{I}_M[O]$ is determined by requiring it to agree
with $O$ at appropriate $M=\dim \tilde{Y}$ interpolation degrees of freedom
$\pi_1, \ldots, \pi_M: Y \to \mathbb{R}$.
I.e., for all $x \in X$ we have

\begin{equation}
	\mathcal{I}_M[O](x) \in \tilde{Y} \quad\text{and}\quad \pi_m(\mathcal{I}_M[O](x)) = \pi_m(O(x)),\ 1\leq m \leq M.
\end{equation}
Due to the locality of finite volume operators, the point evaluations
$\pi_m(O(x))$ can be computed quickly and independently from the dimension of $V_h$.

Since the potential part of $A_\MU$ vanishes identically for solutions
of \eqref{eq:ap4:detailed}, a direct application of empirical interpolation to $O=A_\MU$
results in an unusable approximation, however.
Instead, we further decompose $A_\MU$ and only use empirical interpolation for
appropriate sub-operators.

In the following, we are interested in the behavior of the model in dependence
on the applied current density.
In this case, with $\MU$ being the applied current density, $A_\MU$ decomposes as

\begin{equation}\label{eq:ap4:operator_decomposition}
A_\MU = A_\MU^{(aff)}  + A^{(bv)} + A^{(\ooc)},
\end{equation}
where $A^{(bv)}$, $A^{(\ooc)}$ are the parameter-independent
nonlinear parts of $A_\MU$ corresponding to the Butler-Volmer interface terms
and the summand in \eqref{eq:phiel} containing $\partial \mu / \partial c_{El}$.
Assuming constant $t_+$, the remainder $A^{(aff)}_\MU$ is affine linear
and decomposes as

\begin{equation}\label{eq:ap4:operator_subdecomposition}
	A_\MU^{(aff)} = A^{(const)} + \MU\cdot A^{(bnd)} + A^{(lin)},
\end{equation}
where $A^{(const)}$ is constant and $A^{(bdn)}$, $A^{(lin)}$ are linear,
non-parametric operators.

Now we apply the \textsc{EI-Greedy} algorithm on the training datasets
$\mathcal{S}_{*} = \{A^{(*)}(c^{(n)}_{\MU_s}, \varphi^{(n)}_{\MU_s}),\ 
A^{(*)}(c^{(n, i)}_{\MU_s}, \varphi^{(n, i)}_{\MU_s})\}$, $* \in \{bv, \ooc\}$,
to obtain empirically interpolated operators $\mathcal{I}_{M^{(*)}}[A^{(*)}] \approx A^{(*)}$,
which give us the approximation

\begin{equation}\label{eq:ap4:operator_approximation}
	A_\MU \approx \tilde{A}_\MU = A_\MU^{(aff)}  + \mathcal{I}_{M^{(bv)}}[A^{(bv)}] + \mathcal{I}_{M^{(\ooc)}}[A^{(\ooc)}].
\end{equation}
Substituting \eqref{eq:ap4:operator_approximation} into \eqref{eq:ap4:reduced}
we arrive at the fully reduced model

\begin{align}
\label{eq:ap4:fully_reduced}
		P_{\tilde{V}} \left\{
                \begin{bmatrix}
			\frac{1}{\Delta t^{(n+1)}}(\tilde{c}_{\MU}^{(n+1)} - \tilde{c}_{\MU}^{(n)}) \\
                         0
                \end{bmatrix} 
                 + 
		 \tilde{A}_\MU
                 \left(\begin{bmatrix}
		        \tilde{c}_{\MU}^{(n+1)} \\
			\tilde{\varphi}_{\MU}^{(n+1)}
                 \end{bmatrix}\right) \right\}
         = 0, \notag\\ (\tilde{c}_{\MU}^{(n+1)}, \tilde{\varphi}_{\MU}^{(n+1)}) \in
	 \tilde{V},
\end{align}
with $\tilde{c}^{(0)}_\MU=P_{\tilde{V}_c}$.
After pre-computation of the matrix representations of the linear (constant) operators
$P_{\tilde{V}} \circ A^{(const)}, P_{\tilde{V}} \circ A^{(bnd)}, P_{\tilde{V}} \circ A^{(lin)}: \tilde{V} \to \tilde{V}$,
as well as the projections from the interpolation spaces for $A^{(bv)}$, $A^{(\ooc)}$ onto $\tilde{V}$,
the solution of \eqref{eq:ap4:fully_reduced} can
be obtained quickly for arbitrary new parameters $\MU$ with an effort that only depends on
$\dim\tilde{V}$, $M^{(bv)}$ and $M^{(\ooc)}$.

In the following experiments (see \secref{sec:results:mor}) we are interested in the cell potential as well as the average
\Li concentration in the electrode as functions of time and the applied delithiation current density
$\MU$.
These quantities are linear functionals $s_{cp}, s_{ac}: V_h \oplus V_h \to \mathbb{R}$, assigning
to a state of the cell the respective quantity of interest.
Due to their linearity, the vector representation for the evaluation of $s_{cp}$, $s_{ac}$ on
$\tilde{V}$ can again be pre-computed, such that for any given solution of \eqref{eq:ap4:fully_reduced},
$s_{cp}(\tilde{c}_{\MU}^{(n)}, \tilde{\varphi}_{\MU}^{(n)})$,
$s_{ac}(\tilde{c}_{\MU}^{(n)}, \tilde{\varphi}_{\MU}^{(n)})$ are quickly obtained with
an effort only depending on $\dim \tilde{V}$.

The model order reduction introduces an additional approximation error between full
order and reduced order model that needs to be accounted for in the simulation
workflow.
As we are not aware of any rigorous error estimates which would provide sufficiently
tight error bounds for the model under consideration, we here consider the following
heuristic a posteriori error estimator \cite{ESTIMATOR}:
In addition to $\tilde{V}$ we construct a second, larger validation space
$\hat{V}=\hat{V}_c \oplus \hat{V}_\varphi \supset \tilde{V}_c \oplus \tilde{V}_\varphi = \tilde{V}$
and extended interpolation bases of dimensions $\hat{M}^{(bv)} > M^{(bv)}$ and $\hat{M}^{(1/c)} > M^{(1/c)}$,
yielding a larger reduced model with solutions $(\hat{c}_\mu^{(n)}, \hat{\varphi}_\mu^{(n)}) \in \hat{V}$.
Under the heuristical assumption that

\begin{equation}\label{eq:ap4:saturation}
    \|\hat{c}_\mu^{(n)} - c_\mu^{(n)}\| \leq \Theta \cdot\|\tilde{c}_\mu^{(n)} - c_\mu^{(n)}\|, \quad
    \|\hat{\varphi}_\mu^{(n)} - \varphi_\mu^{(n)}\| \leq \Theta \cdot\|\tilde{\varphi}_\mu^{(n)} - \varphi_\mu^{(n)}\|
\end{equation}
for all timesteps $n$, $\mu \in \mathcal{P}$ with a fixed $\Theta \in [0, 1)$,
from the triangle inequality we immediately obtain the error estimates

\begin{equation}\label{eq:ap4:estimate}
    \|\tilde{c}_\mu^{(n)} - c_\mu^{(n)}\| \leq \frac{1}{1-\Theta} \cdot\|\tilde{c}_\mu^{(n)} - \hat{c}_\mu^{(n)}\|, \quad
    \|\tilde{\varphi}_\mu^{(n)} - \varphi_\mu^{(n)}\| \leq \frac{1}{1-\Theta} \cdot\|\tilde{\varphi}_\mu^{(n)} - \hat{\varphi}_\mu^{(n)}\|.
\end{equation}
The right-hand sides of \eqref{eq:ap4:estimate} can be quickly computed at the
expense of an additional solution of a second (slightly larger) reduced order model.
Note that \eqref{eq:ap4:saturation} is precisely the \emph{saturation assumption} in
the context of hierarchical a posteriori estimates for finite element schemes (see e.g. \cite{BankSmith93}).

\subsection{Algorithmical integration and software interfaces}
\label{sec:workflow:integration}

\MULTIBAT aims to allow computationally fast studies of local effects in the complex microstructure of battery anodes within one software workflow. This is achieved by breaking the multi-disciplinary goal into task units and interfacing these units with \BEST to varying degrees of depth.

\begin{figure}[!ht]
	\centering
	\includegraphics[width=0.8\textwidth]{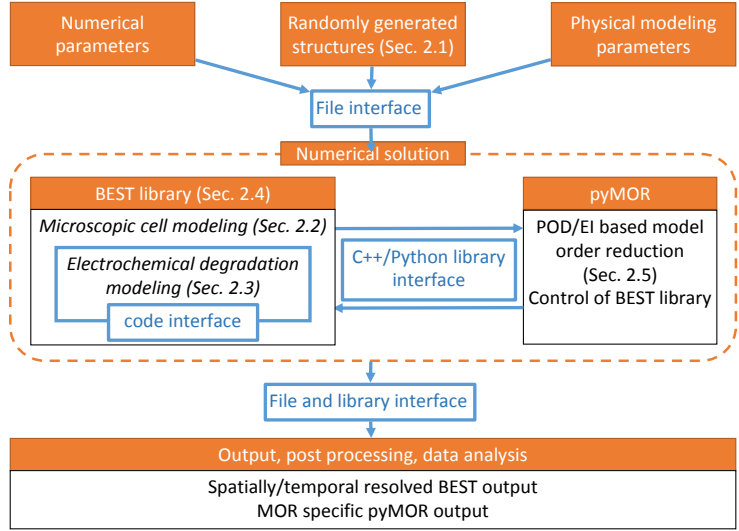}
    \caption{Implementation of the MULTIBAT workflow.}
	\label{fig:workflow}
\end{figure}

The presented workflow resulting from these interfaces is schematically depicted in \figref{fig:workflow}. It allows speeding up the numerical solution of the microscopic cell degradation modeling from \secref{sec:workflow:cell_modeling} and \ref{sec:workflow:degradation_modeling} with discretization from \secref{sec:workflow:discretization} by \BEST, using the randomly generated structures from \secref{sec:workflow:generation} through the \EI based model order reduction approach from \secref{sec:workflow:MOR}. The workflow has been used to create the results and speedups depicted in \secref{sec:results}.

We introduce three distinct interfaces. The first is file-based and allows usage of randomly generated structures of \secref{sec:workflow:generation} in \BEST through a conversion tool. The conversion tool provides standard \BEST geometry input which is matched with physical modeling input parameters and numerical solution parameters suitable for the models of \secref{sec:workflow:cell_modeling} with the extensions from  \secref{sec:workflow:degradation_modeling}.

The second interface allows to extend the \BEST numerical solution code to advanced interface flux modeling between the plated anode-\Li and the electrolyte from \secref{sec:workflow:degradation_modeling}. The model extensions are compiled into the \BEST library.

The third and most extensive interface is library-based and gives \PYMOR runtime access to the solution process, vectors, discretization matrices, Jacobians, linear algebra solver and parameters of \BEST through the \BEST library to carry out \EI based \MOR. The separation is strict: All \MOR-related operations and the Newton methods are carried out in \PYMOR and all evaluations of nonlinear operators and Jacobians are carried out by the \BEST library ordered by \PYMOR.

\section{Workflow demonstration on lithium stripping case study}\label{sec:results}

\subsection{Microstructure generation and electrochemical verification}
\begin{figure}[!ht]
	\centering
		\subfigure[]{
		\includegraphics[width=0.45\textwidth]{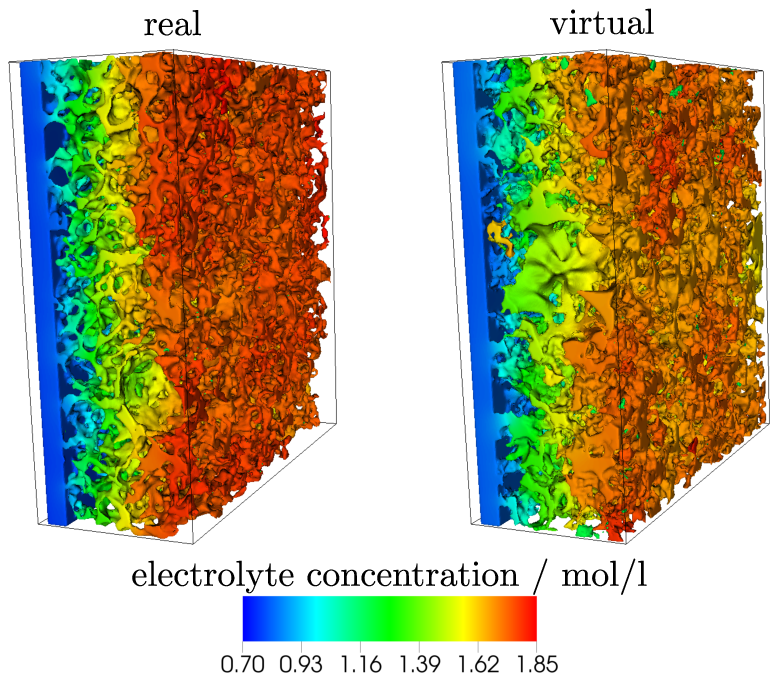}
		\label{fig:res:con:3d}
	}
		\centering
		\subfigure[]{
		\includegraphics[width=0.45\textwidth]{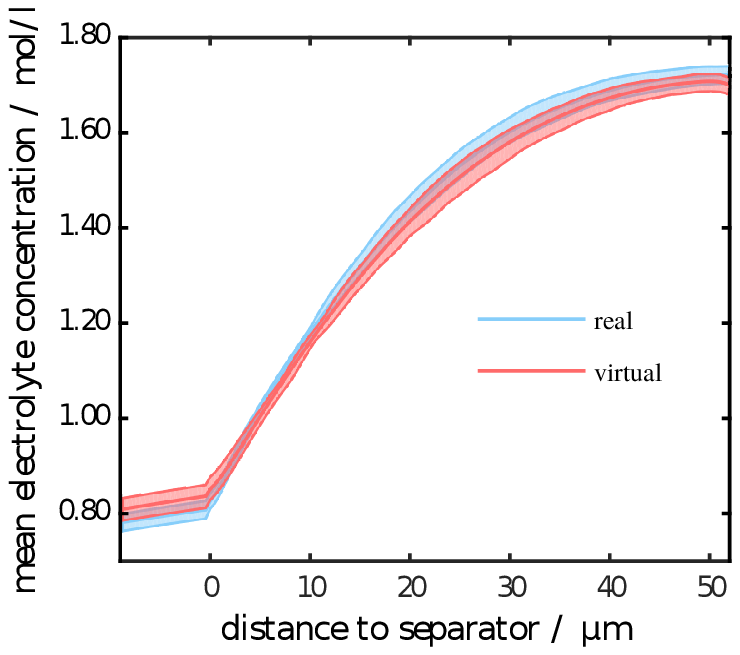}
		\label{fig:res:con:slice}
	}
    \caption{(a) Spatial distribution of electrolyte concentration of a real (left) and virtual (right) microstructure. The same color scale is used (shown below the cutouts). Both structures exhibit larger particles visible as void spaces. Also both cutouts show electrolyte pores, which are less connected to the main pore space: (virtual) Orange part close to the blue and (real) dark red at the upper corner. (b) Mean \Li concentration in the electrolyte as a function of the distance to the separator averaged over the different microstructures. The color shaded areas indicate the $5\%$ and $95\%$-quantiles.  A good accordance can be observed. Reprinted from \cite{Hein2016b} with permission from Elsevier.}
\end{figure}

The stochastic microstructure model used in this work was parameterized on real tomographic image data \cite{Feinauer2015}. As mentioned in \secref{sec:workflow:generation}, the validity of the structural parameterization was investigated through spatially resolved electrochemical simulations \cite{Hein2016b}. The validated stochastic microstructure model is used in this work. In the following a short summary of the electrochemical validation is given. 20 simulated realizations of the stochastic microstructure model and 20 microstructure cutouts from the tomographic image data are used as electrode structure samples for electrochemical simulations. These microstructures are delithiated with a constant current. The simulation results were compared using various electrochemical quantities, such as local current density and \Li concentration. A very good agreement between the real and virtual microstructures was found, see Fig. \ref{fig:res:con:3d}. The advantage of spatially resolved electrochemical simulations is the access to localized inhomogeneities. The spatial distribution of the electrolyte concentration for two cutouts of real and virtual microstructures is shown as an example in \figref{fig:res:con:3d}.

Both cutouts exhibit similar features: less-than-average connected pores and large particles. Apart from the visual similarity between the real and virtual cutouts, averaged quantities were used for a more quantitative comparison. The average \Li concentration in the electrolyte in through direction (from one current collector towards the other) is shown in \figref{fig:res:con:slice}.

The general shapes of the concentration functions are nearly identical. The superposition of transport within the electrolyte and deintercalation of \Li from the solid phase results in a nonlinear gradient. Without any sources of \Li a linear concentration gradient forms in the separator.
More details regarding the electrochemical validation can be found in the corresponding publication \cite{Hein2016b}.

\subsection{Model order reduction}\label{sec:results:mor}
As a first numerical test for the entire developed modeling and simulation workflow (see \secref{sec:workflow}), we simulated the full model (including plated \Li) on a randomly generated half-cell geometry of size $44\si{\micro\metre} \times 44\si{\micro\metre} \times 65.6\si{\micro\metre}$,
which is meshed with a grid of $100 \times 100 \times 149$ voxels (see Fig.~\ref{fig:microsrtucture}).
This size is required to cover a representative volume containing several particles in each direction, and, at the same time, to achieve a sufficient resolution to resolve a relevant part of the electrode's morphology.
Starting with plated \Li and a high \Li concentration in the electrode we performed a delithiation simulation and hence expect to see \Li stripping.
We simulated 60 seconds with constant current densities in the interval
$\mathcal{P}= [2.5, 250]\si{\ampere/\metre^2}$, which corresponds to currents from $\SI{4.84}{\nano\ampere}$ to
$\SI{484}{\nano\ampere}$ or to C-rates from $\mathrm{C}/10$ to $10\mathrm{C}$.

A single simulation of the full order model \eqref{eq:ap4:detailed} requires around \MORDetailedTimeHours~hours
(cf.~Tab.~\ref{tab:mor_times}).
To generate the snapshot data for the computation of the reduced order model \eqref{eq:ap4:fully_reduced},
the full order model \eqref{eq:ap4:detailed} was solved for the three delithiation current densities
$\min \mathcal{P}$, $\max \mathcal{P}$ and $(\min \mathcal{P} + \max \mathcal{P})/2$.
The reduced spaces $\tilde{V}_c$, $\tilde{V}_\varphi$, $\hat{V}_c$, $\hat{V}_\varphi$, as well as the interpolation spaces for
$\mathcal{I}_{M^{(bv)}}[A^{(bv)}]$, $\mathcal{I}_{M^{(\ooc)}}[A^{(\ooc)}]$ were computed using
the POD and \textsc{EI-Greedy} algorithms.
To ensure that a numerically stable reduced order model is obtained, a small relative error tolerance of $10^{-7}$
was chosen, using in each case 97\% of the resulting basis vectors for the construction of
the reduced order model and all basis vectors for construction of the validation model used
for the error estimator \eqref{eq:ap4:estimate}.
The resulting spaces are of the following dimensions:
$\dim \tilde{V}_c = \MORREALCDIM$, $\dim \hat{V}_c = \MORCDIM$,
$\dim \tilde{V}_\varphi = \MORREALPDIM$, $\dim \hat{V}_\varphi = \MORPDIM$,
$M^{(bv)} = \MORREALBVDIM$, $\hat{M}^{(bv)} = \MORBVDIM$,
$M^{(\ooc)} = \MORREALLAMBDADIM$ and $\hat{M}^{(\ooc)} = \MORLAMBDADIM$.

\begin{figure}
	\centering
        \includegraphics{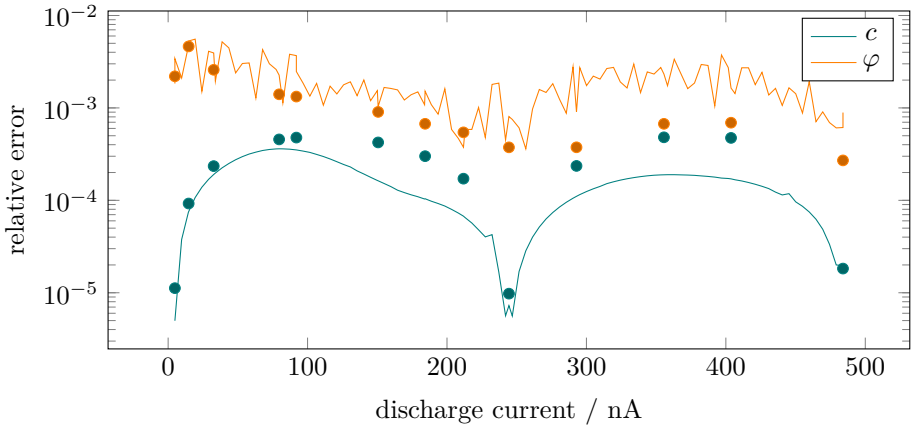}
    \caption{State space model reduction errors (markers) and estimated model reduction errors (solid lines)
    for different applied discharge current densities for the model reduction experiment.
    Plotted is the relative $L^\infty$-in-time $L^2$-in-space error in the concentration ($c$)
	and potential ($\varphi$) variables over a test set of 10 randomly chosen current densities
	$\MU_i \in \mathcal{P} = [2.5, 250]~\si{\ampere/\metre^2}\,, i = 1, \ldots, 10$ in addition to the three
    current densities used for training of the reduced model.
    The model reduction error was estimated for additional 100 equidistantly sampled current densities
    in $\mathcal{P}$, $\Theta = 0$.
    }\label{fig:mor_errors}
\end{figure}

To validate the resulting reduced order model \eqref{eq:ap4:fully_reduced}, we compared the
solutions of \eqref{eq:ap4:fully_reduced} to the full order model \eqref{eq:ap4:detailed} for
10 random parameter values $\MU_i \in \mathcal{P}\,, i=1, \ldots, 10$ in addition to the three snapshot
parameters used for training.
While achieving a relative model reduction error of at most \MORCRelErr~ resp. \MORPRelErr~
for the concentration and potential variables (Fig.~\ref{fig:mor_errors}),
the reduced order model can be simulated in less than \MORReducedTimeMinutes~minutes, yielding
a speedup factor of \MORSpeedup.
Since the generation of the reduced model from the high-dimensional snapshot data is 
faster than a single solution of the full order model, an overall saving of computation time is already
achieved for one additional model simulation (Tab.~\ref{tab:mor_times}).
For $\Theta = 0$ the error estimator \eqref{eq:ap4:estimate} overestimates the real model reduction error in
these 13 parameters by a factor of at most \MORCEffMinRec~(\MORPEffMinRec) for the concentration (potential)
and underestimates the error by a factor of at most \MORCEffMax~(\MORPEffMax).
The estimator was evaluated for 100 additional current densities in $\mathcal{P}$,
yielding a maximum estimated relative error of \MORCEstMax~resp. \MORPEstMax~for concentration and potential.

\begin{figure}
\subfigure[]{
    \centering
        \includegraphics{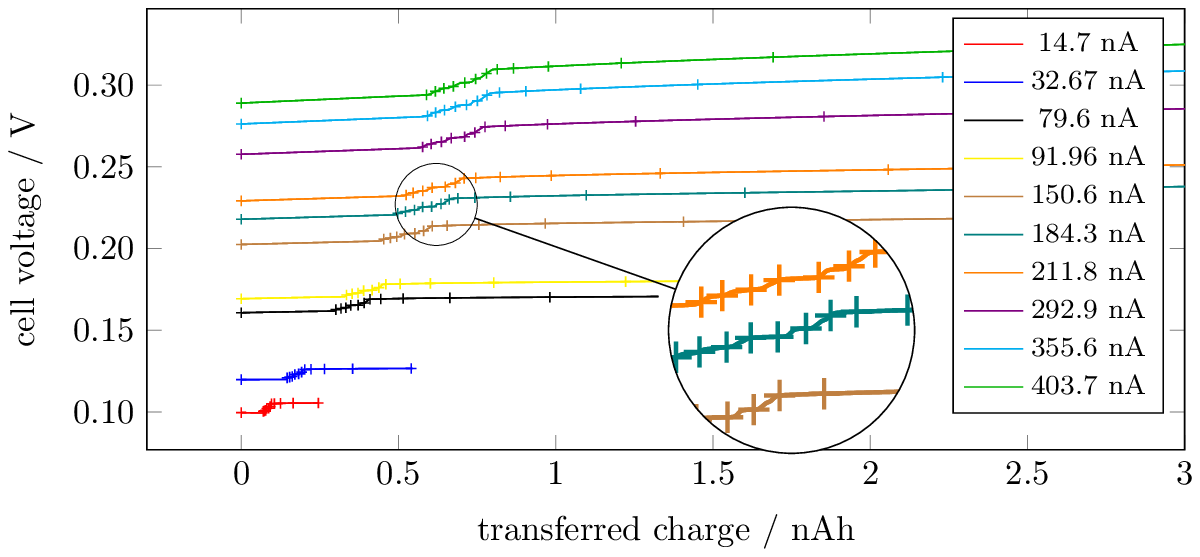}
		\label{fig:ap4:cell_potential}
}
\subfigure[]{
    \centering
        \includegraphics{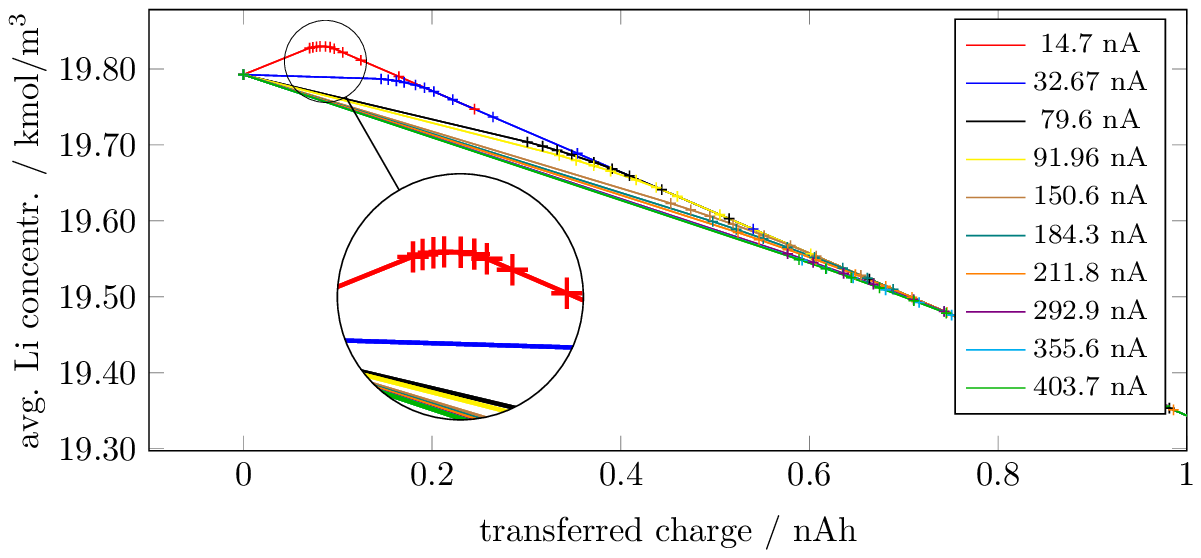}
		\label{fig:ap4:average_concentration}
}
	\caption{(a) Cell voltage over transferred charge for the model reduction experiment 
			for 10 randomly selected current densities $\MU_i$ (cf.\ Fig.~\ref{fig:mor_errors}).
			Solid lines: full model simulation, markers: reduced model simulation (every fifth time step
			marked). (b) Mean \Li concentration inside the solid phase plotted over transferred charge for the model reduction experiment
			 for 10 randomly selected currents densities $\MU_i$ (cf.\ Fig.~\ref{fig:mor_errors}).
			Solid lines: full model simulation, markers: reduced model simulation (every fifth time step
			marked).}
\end{figure}

In Fig.~\ref{fig:ap4:cell_potential} and \ref{fig:ap4:average_concentration}, the cell potential
and average \Li concentration in the electrode have been plotted over the transferred charge
for the 10 random test parameters. A short interpretation of these results is given in the subsequent section. 
Overall, no visual distinction between the data generated by the reduced and full order models can be made.

\begin{table*}
	\caption{Extrapolated timings for the model reduction experiment.
	Time for single full model simulation: \MORDetailedTime\ (median), time for single reduced simulation:
	\MORReducedTime\ (median), time for generation of reduced model from snapshot data: \MOROfflineTime.
	`without \MOR' is the required time if all simulations are performed with the full order model
	\eqref{eq:ap4:detailed}, `with \MOR' is the required time if the reduced order model \eqref{eq:ap4:fully_reduced}
	is used for all simulations after the first 3 snapshot computations (including reduced order model construction).
	All computations have been performed on a single core of an Intel Xeon E5-2698 v3 CPU.}\label{tab:mor_times}
	\begin{center}
		\setlength{\tabcolsep}{0.4ex}
		\newlength{\colsepb}
		\setlength{\colsepb}{8ex}
		\begin{tabular}{c@{\hskip \colsepb}rrr@{\hskip \colsepb}rrr@{\hskip \colsepb}c}\toprule
			simulations & \multicolumn{3}{c@{\hskip \colsepb}}{\makebox[0pt]{without \MOR}} &\multicolumn{3}{c@{\hskip
			\colsepb}}{\makebox[0pt]{\hspace{0.0cm}with \MOR}}  & speedup \\ \midrule
            1&&15h &38m &\multicolumn{3}{c@{\hskip \colsepb}}{--}&--\\
            2&1d &7h &17m &\multicolumn{3}{c@{\hskip \colsepb}}{--}&--\\
            3&1d &22h &55m &\multicolumn{3}{c@{\hskip \colsepb}}{--}&--\\
            4&2d &14h &34m &2d &12h &46m &1.0\\
            5&3d &6h &12m &2d &12h &54m &1.3\\
            10&6d &12h &25m &2d &13h &33m &2.5\\
            50&32d &14h &9m &2d &18h &45m &11.7\\
            100&65d &4h &18m &3d &1h &15m &21.4\\ \midrule
			\multicolumn{7}{l}{\hspace*{0.5cm}limit (= full model vs.\ reduced order model)}& \MORSpeedupB\\	\bottomrule
		\end{tabular}
	\end{center}
\end{table*}

\subsection{The lithium stripping process}
In Fig.~\ref{fig:ap4:cell_potential} the cell voltage is shown for 10 of the applied currents. The cell voltages for all applied currents exhibit a similar shape. A voltage plateau at the start of delithiation is followed by a rise, which is in turn succeeded by a region following the shape of the open-circuit potential $\vart{U}{0}{Graphite}$. The initial voltage plateau results from the stripping reaction (see Eq.~\eqref{eq:ap2:theo:reaction:plating}). The increase in cell voltage begins as soon as the majority of the plated \Li is consumed. The apparent plateau afterwards is the cell voltage of the supporting graphite at about $75\%$ state of charge. This equilibrium potential is shifted with an overpotential, which depends on the applied current. A transferred charge of $\SI{2}{\nano\ampere\hour}$ corresponds to a variation in the state of charge of $4\si{\%}$ since the used microstructure has a maximum capacity of $\SI{49}{nAh}$. Large applied stripping currents lead to a fast decrease of lithium concentration at the surface of graphite. The solid diffusion can not equilibrate the lithium concentration in the electrode in the same rate for large currents as for small currents. 
Therefore, the state of charge at the surface will vary more than the overall change in state of charge. The open-circuit potential will then increase faster for larger currents. Thus resulting in a more sloped cell voltage for the larger currents.
The length of the stripping plateau in the cell voltage depends on the applied current. For smaller currents, the change from constant potential to graphite dominated region is at lower transferred charge. The intercalation of \Li during the stripping of the plated \Li leads to an increase of the \Li concentration in the solid phase, as can be seen in Fig.~\ref{fig:ap4:average_concentration}. For low applied currents a net intercalation during the \Li stripping is visible. As soon as the majority of the plated \Li is dissolved a net delithiation exists. More information about the distribution of the delithiation current on the stripping and intercalation reaction are provided in a recent publication \cite{Hein2016}.

The cell voltages and average \Li concentrations obtained from the model reduction experiments are identical to the ones obtained from the full order model. This indicates that the reduced model sufficiently represents the electrochemical relevant regions in the simulation domain.

\section{Conclusion \& outlook}\label{sec:summary}
We conclude that it is absolutely possible to do a lot of in-depth research on \Liion cells virtually. In this work we have shown one approach to solve many of the existing problems using simulation techniques for the investigation in lithium-ion battery cells.

First, the limitation of 1D or pseudo 2D models which consider only averaged structural quantities and thus neglect all local effects can be overcome by switching to spatially resolved 3D models. Hence we presented a physics-based model that describes the cell's behavior on a microscopic scale and includes effects of lithium plating and stripping. Based on the software tool BEST the mathematical model was implemented and solved in a three-dimensional geometry.

Another current limitation is that the acquisition of tomographic 3D images as basis for simulations is costly and restricts the ability to simulate new structures that have not been produced experimentally. This limitation is overcome by the usage of a 3D stochastic microstructure model which has been implemented in the software library GEOSTOCH. Once the model is fitted to a material by the usage of tomographic 3D images it is possible to generate arbitrary many virtual cutouts with arbitrary sizes. By reasonable changes on the model parameters it is even possible to generate structures that have not been processed experimentally in order to investigate their properties and to test their performance.

Finally, a great problem of simulation-based parameter studies in particular including the plating/stripping behavior is the extensive simulation runtime.
This problem is solved using model order reduction methods implemented in the software library \PYMOR, which speed up the simulation of similar cycles significantly. Tab. \ref{tab:mor_times} shows a speedup of factor \MORSpeedup\ for the simulation of the reduced model in comparison to a full order model simulation.

Overall we have shown that the combinations of all the methods described above work well in a demo scenario and can improve the accuracy of the geometry models, increase the computational speed considerably, and extended the predictive power of the electrochemical battery models.

\section*{Acknowledgement}
This work was partially funded by BMBF under grant numbers 05M13VUA, 05M13PMA, 05M13AMF and 05M13CLA in the programme ``Mathematik f\"{u}r Innovationen in Industrie und Dienstleistungen''.


\begin{thebibliography}{10}

\bibitem{akolkar2013mathematical}
Rohan Akolkar.
\newblock Mathematical model of the dendritic growth during lithium
  electrodeposition.
\newblock {\em Journal of Power Sources}, 232:23--28, 2013.

\bibitem{Arora1999}
P.~Arora, M.~Doyle, and R.~E. White.
\newblock Mathematical modeling of the lithium deposition overcharge reaction
  in lithium-ion batteries using carbon-based negative electrodes.
\newblock {\em Journal of The Electrochemical Society}, 146(10):3543--3553,
  1999.

\bibitem{Arora1998}
P.~Arora, R.~E. White, and M.~Doyle.
\newblock Capacity fade mechanisms and side reactions in lithium-ion batteries.
\newblock {\em Journal of The Electrochemical Society}, 145(10):3647--3667,
  1998.

\bibitem{BankSmith93}
Randolph~E. Bank and R.~Kent Smith.
\newblock A posteriori error estimates based on hierarchical bases.
\newblock {\em SIAM Journal on Numerical Analysis}, 30(4):921--935, 1993.

\bibitem{BarraultMadayEtAl2004}
Maxime Barrault, Yvon Maday, Ngoc~Cuong Nguyen, and Anthony~T. Patera.
\newblock An `empirical interpolation' method: application to efficient
  reduced-basis discretization of partial differential equations.
\newblock {\em Comptes Rendus Mathematique}, 339(9):667--672, 2004.

\bibitem{CaiWhi:2009}
L.~Cai and R.E. White.
\newblock Reduction of model order based on proper orthogonal decomposition for
  lithium-ion battery simulations.
\newblock {\em Journal of The Electrochemical Society}, 156(3):154--161, 2009.

\bibitem{Chiu2013}
S.~N. Chiu, D.~Stoyan, W.~S. Kendall, and J.~Mecke.
\newblock {\em Stochastic Geometry and its Applications}.
\newblock J. Wiley \& Sons, 3rd edition edition, 2013.

\bibitem{Dougherty1992}
E.~Dougherty, editor.
\newblock {\em Mathematical Morphology in Image Processing}.
\newblock Optical Science and Engineering. Taylor \& Francis, 1992.

\bibitem{Doyle1993}
Marc Doyle, Thomas~F Fuller, and John Newman.
\newblock {M}odeling of galvanostatic charge and discharge of the
  lithium/polymer/insertion cell.
\newblock {\em Journal of The Electrochemical Society}, 140(6):1526--1533,
  1993.

\bibitem{Doyle1996}
Marc Doyle, John Newman, Antoni~S. Gozdz, Caroline~N. Schmutz, and Jean-Marie
  Tarascon.
\newblock Comparison of modeling predictions with experimental data from
  plastic lithium ion cells.
\newblock {\em Journal of The Electrochemical Society}, 143(6):1890--1903,
  1996.

\bibitem{DrohmannHaasdonkEtAl2012}
M.~Drohmann, B.~Haasdonk, and M.~Ohlberger.
\newblock Reduced basis approximation for nonlinear parametrized evolution
  equations based on empirical operator interpolation.
\newblock {\em SIAM Journal on Scientific Computing}, 34(2):937--969, 2012.

\bibitem{Ecker2015b}
M.~Ecker, S.~Kabitz, I.~Laresgoiti, and D.~U. Sauer.
\newblock Parameterization of a physico-chemical model of a lithium-ion
  battery: {II}. model validation.
\newblock {\em Journal of The Electrochemical Society}, 162(9):A1849--A1857,
  2015.

\bibitem{Feinauer2015}
Julian Feinauer, Tim Brereton, Aaron Spettl, Matthias Weber, Ingo Manke, and
  Volker Schmidt.
\newblock {Stochastic 3D modeling of the microstructure of lithium-ion battery
  anodes via Gaussian random fields on the sphere}.
\newblock {\em Computational Material Science}, 109:137--146, 2015.

\bibitem{Finegan2016}
Donal~P. Finegan, Mario Scheel, James~B. Robinson, Bernhard Tjaden, Marco~Di
  Michiel, Gareth Hinds, Dan J.~L. Brett, and Paul~R. Shearing.
\newblock Investigating lithium-ion battery materials during overcharge-induced
  thermal runaway: an operando and multi-scale {X}-ray {CT} study.
\newblock {\em Physical Chemistry Chemical Physics}, 18:30912--30919, 2016.

\bibitem{ITWM2014}
{Fraunhofer ITWM}.
\newblock {BEST -- Battery and Electrochemistry Simulation Tool}, 2014.

\bibitem{Fuller1994}
Thomas~F Fuller, Marc Doyle, and John Newman.
\newblock {S}imulation and optimization of the dual lithium ion insertion cell.
\newblock {\em Journal of The Electrochemical Society}, 141(1):1--10, 1994.

\bibitem{Gaiselmann2013}
G.~Gaiselmann, M.~Neumann, L.~Holzer, T.~Hocker, M.~R. Prestat, and V.~Schmidt.
\newblock Stochastic 3{D} modeling of {LSC} cathodes based on structural
  segmentation of {FIB-SEM} images.
\newblock {\em Computational Materials Science}, 67:48--62, 2013.

\bibitem{Gaiselmann2012}
G.~Gaiselmann, R.~Thiedmann, I.~Manke, W.~Lehnert, and V.~Schmidt.
\newblock Stochastic 3{D} modeling of fiber-based materials.
\newblock {\em Computational Materials Science}, 59:75--86, 2012.

\bibitem{Ha14}
Bernard Haasdonk.
\newblock Reduced basis methods for parametrized pdes—a tutorial introduction
  for stationary and instationary problems.
\newblock In Peter Benner, Mario Ohlberger, Albert Cohen, and Karen Willcox,
  editors, {\em Model Reduction and Approximation}, chapter~2, pages 65--136.
  Society for Industrial and Applied Mathematics, 2017.

\bibitem{ESTIMATOR}
Stefan Hain, Mario Ohlberger, Mladjan Radic, and Karsten Urban.
\newblock A hierarchical a-posteriori error estimator for the reduced basis
  method.
\newblock Working paper (under preparation), 2017+.

\bibitem{Hein2016b}
Simon Hein, Julian Feinauer, Daniel Westhoff, Ingo Manke, Volker Schmidt, and
  Arnulf Latz.
\newblock {Stochastic microstructure modeling and electrochemical simulation of
  lithium-ion cell anodes in 3D}.
\newblock {\em Journal of Power Sources}, 336:161--171, dec 2016.

\bibitem{Hein2016}
Simon Hein and Arnulf Latz.
\newblock {Influence of local lithium metal deposition in 3D microstructures on
  local and global behavior of Lithium-ion batteries}.
\newblock {\em Electrochimica Acta}, 201:354--365, 2016.

\bibitem{HesthavenRozzaEtAl2016}
Jan~S Hesthaven, Gianluigi Rozza, and Benjamin Stamm.
\newblock {\em Certified Reduced Basis Methods for Parametrized Partial
  Differential Equations}.
\newblock SpringerBriefs in Mathematics. Springer International Publishing,
  2016.

\bibitem{Hofmann2016}
Tobias Hofmann, Ralf M{\"{u}}ller, Heiko Andr{\"{a}}, and Jochen Zausch.
\newblock {Numerical simulation of phase separation in cathode materials of
  lithium ion batteries}.
\newblock {\em International Journal of Solids and Structures},
  100-101:456--469, 2016.

\bibitem{Hutzenlaub2014}
T.~Hutzenlaub, S.~Thiele, N.~Paust, Robert~M Spotnitz, R.~Zengerle, and
  C.~Walchshofer.
\newblock {Three-dimensional electrochemical Li-ion battery modelling featuring
  a focused ion-beam/scanning electron microscopy based three-phase
  reconstruction of a LiCoO2 cathode}.
\newblock {\em Electrochimica Acta}, 115:131--139, jan 2014.

\bibitem{multibat}
O.~Iliev, A.~Latz, M.~Ohlberger, S.~Schmidt, and V.~Schmidt.
\newblock The multibat project, 2016.

\bibitem{IlievLatzEtAl2012}
Oleg Iliev, Arnulf Latz, Jochen Zausch, and Shiquan Zhang.
\newblock {On some model reduction approaches for simulations of processes in
  Li-ion battery.}
\newblock In {\em {Proceedings of Algoritmy 2012, conference on scientific
  computing, Vysok\'e Tatry, Podbansk\'e, Slovakia}}, pages 161--171. Slovak
  University of Technology in Bratislava, 2012.

\bibitem{Lang2014}
A.~Lang and C.~Schwab.
\newblock Isotropic {G}aussian random fields on the sphere: regularity, fast
  simulation, and stochastic partial differential equations.
\newblock {\em Annals of Applied Probability}, 25(6):3047--3094, 2015.

\bibitem{LV15}
Oliver Lass and Stefan Volkwein.
\newblock Parameter identification for nonlinear elliptic-parabolic systems
  with application in lithium-ion battery modeling.
\newblock {\em Computational Optimization and Applications}, 62(1):217--239,
  2015.

\bibitem{Latz2011a}
A.~Latz and J.~Zausch.
\newblock Thermodynamic consistent transport theory of {L}i-ion batteries.
\newblock {\em Journal of Power Sources}, 196(6):3296--3302, 2011.

\bibitem{Latz2011b}
A.~Latz, J.~Zausch, and O.~Iliev.
\newblock Modeling of species and charge transport in {L}i--ion batteries based
  on non-equilibrium thermodynamics.
\newblock In I.~Dimov, S.~Dimova, and N.~Kolkovska, editors, {\em Numerical
  Methods and Applications}, volume 6046 of {\em Lecture Notes in Computer
  Science}, pages 329--337. Springer, 2011.

\bibitem{Latz2013}
Arnulf Latz and Jochen Zausch.
\newblock {Thermodynamic derivation of a Butler–Volmer model for
  intercalation in Li-ion batteries}.
\newblock {\em Electrochim. Acta}, 110:358--362, 2013.

\bibitem{Latz2015}
Arnulf Latz and Jochen Zausch.
\newblock {Multiscale modeling of lithium ion batteries: Thermal aspects}.
\newblock {\em Beilstein J. Nanotechnol.}, 6:987--1007, 2015.

\bibitem{Lautensack2007}
C.~Lautensack.
\newblock {\em Random {L}aguerre Tessellations}.
\newblock PhD thesis, Universität Karlsruhe (TH), 2007.

\bibitem{Legrand2014a}
N.~Legrand, B.~Knosp, P.~Desprez, F.~Lapicque, and S.~Ra{\"{e}}l.
\newblock {Physical characterization of the charging process of a Li-ion
  battery and prediction of Li plating by electrochemical modelling}.
\newblock {\em Journal of Power Sources}, 245:208--216, jan 2014.

\bibitem{Less2012a}
G.~B. Less, J.~H. Seo, S.~Han, Ann~Marie Sastry, Jochen Zausch, Arnulf Latz,
  Sebastian Schmidt, C.~Wieser, D.~Kehrwald, and S.~Fell.
\newblock Micro-scale modeling of {L}i-ion batteries: Parameterization and
  validation.
\newblock {\em Journal of The Electrochemical Society}, 159(6):697--704, 2012.

\bibitem{Cai2009}
C.~Long and Ralph~E. White.
\newblock Reduction of model order based on proper orthogonal decomposition for
  lithium-ion battery simulations.
\newblock {\em Journal of The Electrochemical Society}, 156(3):154--161, 2009.

\bibitem{Mayer2004}
J.~Mayer, V.~Schmidt, and F.~Schweiggert.
\newblock A unified simulation framework for spatial stochastic models.
\newblock {\em Simulation Modelling Practice and Theory}, 12:307--326, 2004.

\bibitem{MRS16}
R.~Milk, S.~Rave, and F.~Schindler.
\newblock {pyMOR} - {G}eneric algorithms and interfaces for model order
  reduction.
\newblock {\em SIAM J. Sci. Comput.}, 35(5):194--216, 2016.

\bibitem{Moller1989}
J.~M{\o}ller.
\newblock Random tessellations in $\mathbb{R}^d$.
\newblock {\em Advances in Applied Probability}, 21(1):37--73, 1989.

\bibitem{monroe2003dendrite}
Charles Monroe and John Newman.
\newblock Dendrite growth in lithium/polymer systems a propagation model for
  liquid electrolytes under galvanostatic conditions.
\newblock {\em Journal of The Electrochemical Society}, 150(10):1377--1384,
  2003.

\bibitem{Neumann2016}
M.~Neumann, J.~Stanek, O.~Pecho, L.~Holzer, V.~Benes, and V.~Schmidt.
\newblock {Stochastic 3D modeling of complex three-phase microstructures in
  SOFC-electrodes with completely connected phases}.
\newblock {\em Computational Materials Science}, 118:353--364, 2016.

\bibitem{Newman2003}
J.~Newman, K.~E. Thomas, H.~Hafezi, and D.~R. Wheeler.
\newblock Modeling of lithium-ion batteries.
\newblock {\em Journal of Power Sources}, 119:838--843, 2003.

\bibitem{Newman2004}
John Newman and Karen~E. Thomas-Alyea.
\newblock {\em {Electrochemical systems}}.
\newblock J. Wiley {\&} Sons, 3rd edition edition, 2004.

\bibitem{Newman2012}
John Newman and Karen~E Thomas-Alyea.
\newblock {\em {E}lectrochemical {S}ystems}.
\newblock John Wiley \& Sons, {T}hird edition, 2012.

\bibitem{Newman1975}
John Newman and William Tiedemann.
\newblock {Porous-electrode theory with battery applications}.
\newblock {\em AIChE Journal}, 21(1):25--41, 1975.

\bibitem{OR16a}
M.~Ohlberger and S.~Rave.
\newblock Localized reduced basis approximation of a nonlinear finite volume
  battery model with resolved electrode geometry.
\newblock In P.~Benner, M.~Ohlberger, A.~Patera, G.~Rozza, and K.~Urban,
  editors, {\em Model Reduction of Parametrized Systems}, number~17 in MS\&A,
  pages 201--212. Springer International Publishing, 2017.

\bibitem{OhlbergerRaveEtAl2016}
M.~Ohlberger, S.~Rave, and F.~Schindler.
\newblock Model reduction for multiscale lithium-ion battery simulation.
\newblock In Bülent Karasözen, Murat Manguoğlu, Münevver Tezer-Sezgin,
  Serdar Göktepe, and Ömür Uğur, editors, {\em Numerical Mathematics and
  Advanced Applications ENUMATH 2015}, volume 112 of {\em Lecture Notes in
  Computational Science and Engineering}, pages 317--331. Springer, 2016.

\bibitem{OhlbergerRaveEtAl2014}
Mario Ohlberger, Stephan Rave, Sebastian Schmidt, and Shiquan Zhang.
\newblock A model reduction framework for efficient simulation of li-ion
  batteries.
\newblock In J.~Fuhrmann, M.~Ohlberger, and C.~Rohde, editors, {\em Finite
  Volumes for Complex Applications VII-Elliptic, Parabolic and Hyperbolic
  Problems}, volume~78 of {\em Springer Proceedings in Mathematics \&
  Statistics}, pages 695--702. Springer International Publishing, 2014.

\bibitem{Popov2011}
P.~Popov, Y.~Vutov, S.~Margenov, and O.~Iliev.
\newblock Finite volume discretization of equations describing nonlinear
  diffusion in {L}i-ion batteries.
\newblock {\em Numerical Methods and Applications}, 6046:338--346, 2011.

\bibitem{Prim1957}
Robert~Clay Prim.
\newblock Shortest connection networks and some generalizations.
\newblock {\em The Bell System Technical Journal}, 36(6):1389--1401, 1957.

\bibitem{PYMOR}
{pyMOR Developers and Contributors}.
\newblock {pyMOR -- Model Order Reduction with Python}, 2013--2016.

\bibitem{QuarteroniManzoniEtAl2016}
Alfio Quarteroni, Andrea Manzoni, and Federico Negri.
\newblock {\em Reduced Basis Methods for Partial Differential Equations}.
\newblock La Matematica per il 3+2. Springer International Publishing, 2016.

\bibitem{Remmlinger2014}
J{\"u}rgen Remmlinger, Simon Tippmann, Michael Buchholz, and Klaus Dietmayer.
\newblock {Low-temperature charging of lithium-ion cells Part II: Model
  reduction and application}.
\newblock {\em Journal of Power Sources}, 254:268--276, 2014.

\bibitem{Santhanagopalan2006}
Shriram Santhanagopalan, Qingzhi Guo, Premanand Ramadass, and Ralph~E White.
\newblock Review of models for predicting the cycling performance of lithium
  ion batteries.
\newblock {\em Journal of Power Sources}, 156(2):620--628, 2006.

\bibitem{sirovich87}
L.~Sirovich.
\newblock Turbulence and the dynamics of coherent structures {P}art {I}:
  {C}oherent structures.
\newblock {\em Quarterly of Applied Mathematics}, 45(3):561--571, 1987.

\bibitem{Smith2009}
Madeleine Smith, R.~Edwin Garc{\'{i}}a, and Quinn~C. Horn.
\newblock The effect of microstructure on the galvanostatic discharge of
  graphite anode electrodes in {LiCoO$_2$}-based rocking-chair rechargeable
  batteries.
\newblock {\em Journal of The Electrochemical Society}, 156(11):A896--A904,
  2009.

\bibitem{SKTOJS_12}
O.~Stenzel, L.J.A. Koster, R.~Thiedmann, S.~D. Oosterhout, R.~A.~J. Janssen,
  and V.~Schmidt.
\newblock A new approach to model-based simulation of disordered polymer blend
  solar cells.
\newblock {\em Advanced Functional Materials}, 22:1236--1244, 2012.

\bibitem{Tang2009}
Maureen Tang, Paul Albertus, and John Newman.
\newblock Two-dimensional modeling of lithium deposition during cell charging.
\newblock {\em Journal of The Electrochemical Society}, 156(5):390--399, 2009.

\bibitem{Taralov2014}
Maxim Taralov, V~Taralova, P~Popov, O~Iliev, Arnulf Latz, and Jochen Zausch.
\newblock On 2{D} finite element simulation of a thermodynamically consistent
  {L}i-ion battery microscale model.
\newblock In Angela Slavova, editor, {\em Mathematics in Industry}, pages
  148--161. Cambridge Scholars Publishing, Newcastle upon Tyne, 2014.

\bibitem{Tippmann2013}
S.~Tippmann, D.~Walper, B.~Spier, and W.~G. Bessler.
\newblock Low-temperature charging of lithium-ion cells. {P}art {I}:
  {E}lectrochemical modeling and experimental investigation on degradation
  behavior.
\newblock {\em Journal of Power Sources}, 252:305--316, 2014.

\bibitem{Torquato2002}
S.~Torquato.
\newblock {\em Random Heterogeneous Materials: Microstructure and Macroscopic
  Properties}.
\newblock Springer, 2002.

\bibitem{Vetter2005}
J.~Vetter, P.~Novak, M.~R. Wagner, C.~Veit, K.-C. M{\"o}ller, J.~O. Besenhard,
  M.~Winter, M.~Wohlfahrt-Mehrens, C.~Vogler, and A.~Hammouche.
\newblock Ageing mechanisms in lithium-ion batteries.
\newblock {\em Journal of Power Sources}, 147(1):269--281, 2005.

\bibitem{Wang2007}
Chia-Wei Wang and Ann~Marie Sastry.
\newblock Mesoscale modeling of a {L}i-ion polymer cell.
\newblock {\em Journal of The Electrochemical Society}, 154(11):1035--1047,
  2007.

\bibitem{VW13}
A.~Wesche and S.~Volkwein.
\newblock The reduced basis method applied to transport equations of a
  lithium-ion battery.
\newblock {\em COMPEL: The International Journal for Computation and
  Mathematics in Electrical and Electronic Engineering}, 32:1760--1772, 2013.

\bibitem{Westhoff2016}
D.~Westhoff, J.~Feinauer, K.~Kuchler, T.~Mitsch, I.~Manke, S.~Hein, A.~Latz,
  and V.~Schmidt.
\newblock Parametric stochastic 3{D} model for the microstructure of anodes in
  lithium-ion power cells.
\newblock {\em Computational Materials Science}, 126:453--467, 2017.

\bibitem{Westhoff2015}
D.~Westhoff, J.~J. van Franeker, T.~Brereton, D.~P. Kroese, R.~A.~J. Janssen,
  and V.~Schmidt.
\newblock Stochastic modeling and predictive simulations for the microstructure
  of organic semiconductor films processed with different spin coating
  velocities.
\newblock {\em Modelling and Simulation in Materials Science and Engineering},
  23:045003, 2015.

\bibitem{Yan2013}
Bo~Yan, Cheolwoong Lim, Leilei Yin, and Likun Zhu.
\newblock {Simulation of heat generation in a reconstructed LiCoO2 cathode
  during galvanostatic discharge}.
\newblock {\em Electrochimica Acta}, 100:171--179, jun 2013.

\end{thebibliography}
\end{document}